\documentclass[11pt]{amsart}

\usepackage{graphicx}
\usepackage{float}
\usepackage{natbib} 
\usepackage{dsfont}

\newtheorem{thm}{Theorem}[section]
\newtheorem{lem}{Lemma}[section]
\newtheorem{prop}{Proposition}[section]
\newtheorem{cor}{Corollary}[section]
\newtheorem{Defi}{Definition}[section]

\def\1{\mathds{1}}
\newcommand{\pr}{\mathbb{P}}
\newcommand{\E}{\mathbb{E}}
\newcommand{\R}{\mathbb{R}}
\def\Var{{\rm Var}}

\begin{document}

\title{Adaptive pointwise estimation for pure jump L\'evy processes}
\author{M\'elina Bec*, Claire Lacour**}
\thanks{*  UMR CNRS 8145 MAP5, Universit\'e Paris Descartes, ** Laboratoire de Math\'ematiques d'Orsay, Universit\'e Paris-Sud}
 
 \begin{abstract}
 
 This paper is concerned with adaptive kernel estimation of the L\'evy density $N(x)$ for bounded-variation pure-jump L\'evy processes. The sample path is observed at $n$ discrete instants in the "high frequency" context ($ \Delta $ =  $ \Delta(n) $ tends to zero while $n \Delta $ tends to infinity). We construct a collection of kernel estimators of the function $g(x)=xN(x)$ and propose a method of local adaptive selection of the bandwidth. 
 We provide an oracle inequality and a rate of convergence for the quadratic pointwise risk. This rate is proved to be the optimal minimax rate. We give examples and simulation results for processes fitting in our framework. We also consider the case of irregular sampling. 
\bigskip


  \end{abstract}

 \maketitle

\begin{center}
 \today
\end{center}
 \bigskip

  \section{Introduction}\label{Introduction}
Consider $(L_t, t \geq 0)$ a real-valued L\'evy process with characteristic function given by:
\begin{equation} \label{fc}
\psi_{t}(u)= {\mathbb E}(\exp{i u L_t})=\exp{(t \int_{{\mathbb R}}(e^{iux} -1) N(x)dx)}.
\end{equation}
We assume that the L\'evy measure admits a density $N$ and that the function $g(x)=xN(x) $ is integrable. 
Under these assumptions, $(L_t, t \geq 0)$ is a pure jump L\'evy process without drift and with finite variation on compact sets. Moreover $\mathbb{E}(|L_t|) < \infty $ (see \cite{Ber}). Suppose that we have discrete observations $(L_{k\Delta},{k=1,...,n})$ with sampling interval $\Delta$. Our aim in this paper is the nonparametric adaptive kernel estimation of the function $g(x)=xN(x)$ based on these observations under the asymptotic framework $n$ tends to $\infty$. This subject has been recently investigated by several authors. 
\cite{FH} 
 use a penalized projection method to estimate the L\'evy density on a compact set separated from $0$. Other authors develop an estimation procedure based on empirical estimations of the characteristic function $\psi_{\Delta}(u) $ of the increments $ (Z_k^{\Delta}= L_{k\Delta}- L_{(k-1)\Delta}, k=1, \ldots, n)$ and its derivatives followed by a Fourier inversion to recover the L\'evy density.
For low frequency data ($\Delta$ is fixed), we can quote 
\cite{WK}
, or \cite{JM} 
 for a parametric study. Still in the low frequency framework,
\cite{NR} 
 estimate $\nu(x)=x^2N(x)$ in the more general case with drift and volatility, 
and \cite{CGC} use model selection to build an adaptive estimator. 
An adaptive method to estimate linear functionals is also given in \cite{Kappus}.
\cite{Belomestny11} addresses the issue of inference for time-changed L\'evy processes with results in term of uniform and pointwise distance.

In the high frequency context, which is our concern in this paper, the problem is simpler since, for any fixed $u$,  $\psi_{\Delta}(u) \rightarrow 1 $ when $\Delta \rightarrow 0 $. This implies that $ \psi_{\Delta}(u)$ need not to be estimated and can simply be replaced by $1$ in the estimation procedures. This is what is done in \cite{CG}. These authors start from the equality:
\begin{equation}\label{equationn1}
\mathbb{E}\left[ Z_k^{\Delta} e^{iuZ_k^{\Delta}} \right] = -i \psi'_{\Delta}(u)=\Delta\psi_{\Delta}(u) g^*(u),
\end{equation}
obtained by differentiating \eqref{fc}.
Here $g^*(u)=\int e^{iux}g(x) dx$ is the Fourier transform of $g$, well defined since we assume $g$ integrable.
Then, as $ \psi_{\Delta}(u) \simeq 1$, equation (\ref{equationn1}) writes $\mathbb{E}\left[ Z_k^{\Delta} e^{iuZ_k^{\Delta}} \right] \simeq \Delta g^*(u)$.
This gives an estimator of $g^*(u) $ as follows:
\begin{equation*}
\frac{1}{n\Delta}\sum_{k=1}^{n} Z_k^{\Delta} e^{iuZ_k^{\Delta}}.
\end{equation*}
Now, to recover $g$, the authors apply Fourier inversion with cutoff parameter $m$. Here, we rather introduce a kernel to make inversion possible:
\begin{equation*}
\frac{1}{n\Delta}\sum_{k=1}^{n} Z_k^{\Delta} K^*(uh) e^{iuZ_k^{\Delta}} \label{estg*}
\end{equation*}
which is in fact the Fourier transform of $1/(nh\Delta) \sum_{k=1}^{n} Z_k^{\Delta} K((x-Z_k^{\Delta}) / h)$.
At the end, in the high frequency context, a direct method without Fourier inversion can be applied. 
Indeed, a consequence of (\ref{equationn1}) is that the empirical distribution:
\begin{equation*}
\hat\mu_n(dz)=\frac{1}{n\Delta} \sum_{k=1}^{n} Z_k^{\Delta} \delta_{Z_k^{\Delta}}(dz)
\end{equation*}
weakly converges to $g(z)dz$ (note that the idea of exploiting this weak convergence is already present in \cite{FigueroaCT09}). This suggests to consider kernel estimators of $g$ of the form
\begin{equation}\label{estimateur}
\hat g_h (x)= K_h \star \hat\mu_n (x) = \frac{1}{n\Delta} \sum_{k=1}^{n} Z_k^{\Delta} K_h( x-Z_k^{\Delta})
\end{equation}
where $K_h(x)=(1/h) K(x/h) $ and
$K$ is a kernel such that $\int K=1$.
Below, we study the quadratic pointwise risk of the estimators $\hat g_h(x)$ and evaluate the rate of convergence of this risk as $n$ tends to infinity, $\Delta=\Delta(n)$ tends to $0$ and $h= h(n)$ tends to 0.
This is done under H\"{o}lder regularity assumptions for the function $g$. Note that a pointwise study involving a kernel estimator can be found in 
\cite{Gugu} 
for more specific compound Poisson processes, but the estimator is different from ours, as well as the observation scheme. 
In \cite{Figueroa2011} a pointwise central limit theorem is given for the estimation of the L\'evy density, as well as confidence intervals.
Still in the high frequency context, we can cite \cite{duval12} for the estimation of a compound Poisson process with low conditions on $\Delta$, but for integrated distance. 

In this paper, we study local adaptive bandwidth selection (which the previous authors do not consider). For a given non-zero real $x_0$, we
select a bandwidth $\hat h(x_0)$ such that the resulting adaptive estimator $ \hat g_{\hat h(x_0)}(x_0)$  automatically reaches the optimal rate of convergence corresponding to the unknown regularity of the function $g$. The method of bandwidth selection follows the scheme   
  developped by \cite{GL11} for density estimation. The advantage of our kernel method is that it allows us to estimate the L\'evy density at a fixed point, with a local adaptive choice. This method is 
  easy to implement, and we show its good numerical performance on different examples.
  Moreover our contribution includes an alternative proof for a lower bound result (see \cite{Figueroa09}) which proves the optimality of the rate for this pointwise estimation.
We also study the framework of irregular sampling. 

In Section~\ref{Estimator}, we give notations and assumptions. In Section~\ref{Risk}, we study the pointwise mean square error (MSE) of $\hat g_h (x_0) $ given in (\ref{estimateur}) for $g$ belonging to a H\"{o}lder class of regularity $\beta$ and we present the bandwidth selection method together with both lower and upper risk bound for our adaptive estimator. The rate of convergence of the risk is $(\log(n\Delta)/n\Delta)^{2\beta/2\beta+1}$ which is expected in adaptive pointwise context. Examples and simulations in our framework are discussed in Section~\ref{Examples}. The case of irregular sampling is addressed in Section~\ref{IS}  and proofs are gathered in Section~\ref{Proofs}.

 \section{Notations and assumptions}\label{Estimator}
 We present the assumptions on the kernel $K$ and on the function $g$ required to study the estimator given by (\ref{estimateur}).
 First, we set some notations.
For any functions $u,v$, we denote by
$u^*$ the Fourier transform of $u$, 
$u^*(y)=\int e^{iyx}u(x)dx$ and by $\|u\|$, $<u,v>$, $u\star v$ the quantities
$$\|u\|^2=\int |u(x)|^2 dx,$$ 
$$ <u,v>= \int u(x)\overline{v}(x)dx\mbox{ with }
z\overline{z}=|z|^2 \mbox{ and } u\star v(x)=\int u(y) v(x-y)dy.$$
For a positive real $\beta$, $\lfloor \beta \rfloor $ denotes the largest integer strictly smaller than $ \beta $.
Let us also define the following functional space:
\begin{Defi} \label{Concent} (H\"{o}lder class)
Let $ \beta> 0 $, $ L> 0 $ and let $ l = \lfloor \beta \rfloor $. 
The H\"{o}lder class $ \mathcal{H} (\beta, L) $ on $ \mathbb {R} $  is the set of all functions  $ f: \mathbb {R} \longrightarrow \mathbb{R} $ such that derivative $ f ^ {(l)} $ exists and verifies:
\begin{eqnarray*} \label{holder}| f^{(l)}(x)-f^{(l)}(y) | \leq L|x-y|^{\beta - l} ,\quad\forall x,y\in \mathbb {R}.
\end{eqnarray*} 
 \end{Defi}

We can now define the assumptions concerning the target function $g$:
 \begin{description}
 \item[G1]  $g\in \mathbb{L}^2 $
 \item[G2] $g^*$ is differentiable almost everywhere and its derivative belongs to $\mathbb{L}^1$
 \item[G3($p$)] For $p$ integer, $\int |x|^{p-1}|g(x)|dx <\infty$
\item[G4($\beta$)] $ g \in \mathcal{H}(\beta,L) $ 
\item[G5] $ g ' $ exists and is uniformly bounded 
\end{description}

The first assumption is natural to use Fourier analysis, as well as  G3($1$). Assumption G3($p$) ensures that $\E|Z_1^{\Delta}|^p<\infty$.
G4 is a classical regularity assumption in nonparametric estimation; it allows to quantify the bias (see \cite{TSY}).
Note that G5 implies that $g\in \mathcal{H}(1, L')$ so we can assume $\beta\geq 1$.

Now let us describe which kind of kernel we choose for our estimator. 
For $ m \geq $ 1 an integer,  we say that $ K: \mathbb{R} \rightarrow \mathbb {R} $ is a kernel of order $m$ if functions 
$ u \mapsto u ^ j K (u), j = 0,1,. .., m$ are integrable and satisfy
\begin{eqnarray} \label{noy1}
&&\int K(u)du= 1, \qquad
\int u^j K(u) du = 0, \quad j \in \lbrace 1,...,m\rbrace  .
\end{eqnarray}
Let us  define the following conditions
\begin{description}
 \item[K1] $K$ belongs to $\mathbb{L}^1\cap\mathbb{L}^2\cap \mathbb{L}^\infty$ and $K^*\in \mathbb L^1$
\item[K2($\beta$)] The kernel 
$K$ is of order $l=\lfloor \beta \rfloor $ and $\int |x|^\beta |K(x)| dx < + \infty$
\end{description}

These assumptions  are standard when working on problems of estimation by kernel methods.
Note that there is a way to build a kernel of order $l$. Indeed, let $u$ be a bounded integrable function such that $u \in \mathbb{L}^2$, $u^* \in \mathbb{L}^1$ and $\int u(y)dy =1$, and set for any given integer $l$,    
\begin{eqnarray}\label{noyauOP}
K(t) = \sum_{k=1}^l \binom{l}{k} (-1)^{k+1} \frac{1}{k} u \left( \frac{t}{k}\right).
\end{eqnarray}
The kernel $K$ defined by (\ref{noyauOP}) is a kernel of order $l$ which also satisfies K1
(see \cite{KLP} and  \cite{GL11}). 
As usual, we define $K_h$ by 
$$\forall x\in \mathbb{R} \qquad K_h(x)=\frac1h K\left(\frac{x}h\right).$$

In all the following we fix $x_0 \in \mathbb{R}$, $x_0\neq 0$.

\section{Risk bound}\label{Risk}

\subsection{Risk bound for a fixed bandwidth}
In this subsection, the bandwidth $h$ is fixed, thus we omit the subscript $h$ for the sake of simplicity: we denote $\hat g=\hat g_h$.
 The usual bias variance decomposition of the Mean Squared Error yields:
\begin{eqnarray*}
MSE(x_0,h):=\mathbb{E}[{(\widehat{g}(x_0)-g(x_0))}^2] = \mathbb{E}[{(\widehat{g}(x_0)-\mathbb{E}[\widehat{g}(x_0)}])^2] + {(\mathbb{E}[\widehat{g}(x_0)]-g(x_0))}^2. 
 \end{eqnarray*}
But the bias needs further decomposition:
\begin{eqnarray*}
b(x_0)^2:={(\mathbb{E}[\widehat{g}(x_0)]-g(x_0))}^2 & \leq &  2{b_1(x_0)}^2 + 2{b_2(x_0)}^2
\end{eqnarray*} 
with the usual bias, 
\begin{eqnarray*}
b_1(x_0)= K_h\star g(x_0)-g(x_0),
\end{eqnarray*}
and the bias resulting from the approximation of $ \psi_{\Delta}(u)$ by 1,
\begin{eqnarray*}
b_2(x_0)= \mathbb{E}[\widehat{g}(x_0)]-K_h\star g(x_0).
\end{eqnarray*} 
We can provide the following bias bound:

\begin{lem} \label{monlem1}
Under G3(1), G4($\beta$), G5 and if the kernel $K$ satisfies K1 and K2($\alpha$) with $\alpha\geq \beta$ 
\begin{eqnarray*}
|b(x_0)|^2 & \leq &  c_1 h^{2\beta} +  {c'}_1\Delta ^2
\end{eqnarray*}
with $c_1=  2 \left( { L}/{\lfloor \beta \rfloor !} \int |K(v)| |v|^\beta dv \right) ^2$ and ${c'}_1= 2 (2\|g'\|_\infty{\| g \|}_1 {\| K \|}_1)^2$.
\end{lem}
Moreover, the variance is controlled as follows:

\begin{lem} \label{monlem2}
Under G1 and G2,  and if the kernel satisfies K1, we have
\begin{eqnarray*}
\Var[\widehat{g}(x_0)] \leq \frac{1}{nh\Delta}  \frac{\| K \|^2_2}{2\pi} (\|(g^*)'\|_1  + {\| g ^*\|}^2_2 \Delta) 
\leq c_2 \frac{1}{nh \Delta} + {c'}_2 \frac{1}{nh}
\end{eqnarray*}
with $c_2 = \|(g^*)'\|_1{\| K \|}^2_2 /(2\pi) $ and ${c'}_2= {\| K \|}^2_2\| g \|^2_2 $.
\end{lem}

Lemmas \ref{monlem1} and \ref{monlem2} lead us to the following risk bound:

\begin{prop} \label{proprisque1}
Under G1, G2, G3(1), G4($\beta$), G5 and if $K$ satifies K1 and K2($\alpha$) with $\alpha\geq \beta$, we have
\begin{eqnarray} \label{borne}
MSE(x_0,h) \leq c_1 h^{2\beta} + c_2 \frac{1}{nh \Delta} + {c'}_2 \frac{1}{nh} + {c'}_1 \Delta^2.
\end{eqnarray}
\end{prop} 

Recall that $\Delta= \Delta(n)$ is such that $\lim_{n\rightarrow +\infty} \Delta = 0$, thus $1/nh$ is negligible compared to $1/nh \Delta$. 
For the two first terms 
the optimal choice of $h$ is $h_{opt} \propto((n\Delta)^{-\frac{1}{2\beta + 1}})$ and the associated rate has order
$O\left((n\Delta)^{-\frac{2\beta}{2\beta + 1}}\right) $. 
Next, a sufficient condition for $\Delta^2\leq (n\Delta)^{-\frac{2\beta}{2\beta + 1}}$ for all $\beta$ is 
\begin{eqnarray} \label{cond}
\Delta=O( n^{-1/3}) .
 \end{eqnarray}
\begin{prop} \label{proprisque2}
Under the assumptions of Proposition \ref{proprisque1} and under condition (\ref{cond}), the choice $h_{opt}\propto((n\Delta)^{-\frac{1}{2\beta + 1}})$  minimizes the risk bound (\ref{borne}) and gives $MSE(x_0,h_{opt})=O((n\Delta)^{-\frac{2 \beta}{2\beta + 1}})$.
As a consequence 
$\mathbb{E}[{(\widehat{g}(x_0)/x_0-N(x_0))}^2] =O((n\Delta)^{-\frac{2 \beta}{2\beta + 1}})$.
\end{prop}
We can link this result to the one of \cite{Figueroa2011} who proves that his projection estimator $\widehat{N}$ is such that 
${(\widehat{N}(x_0)-N(x_0))}(n\Delta)^{\alpha}$ tends to a normal distribution for any $0<\alpha<{ \beta}/(2\beta + 1)$. 


%

The rate obtained in Proposition~\ref{proprisque2} turns out to be the optimal minimax rate of convergence over the class  ${\mathcal H}(\beta,L)$. 
This result is proved in \cite{Figueroa09} in the more general case of estimators based on the whole path of the process up to time $n\Delta$. In our case of discrete sampling, another proof is given in Section~\ref{preuveborneinf}, where we prove the following result:
 
 \begin{thm}\label{lower} 
 Assume $\Delta=O(1)$ and $\Delta^{-1}=O(n)$. Let $x_0\neq 0$. 
 There exists $C>0$ such that for any estimator $\hat g_n(x_0)$ based on observations $Z_1^{\Delta},\dots, Z_n^{\Delta}$, and for $n$ large enough,
 $$\sup_{g\in {\mathcal H}(\beta,L)}{\mathbb E}_g\left[(\hat g_n(x_0)-g(x_0))^2\right] \geq C (n\Delta)^{-\frac{2\beta}{2\beta+1}}.$$
 \end{thm}
 Obviously, the result is also true replacing $g$ by the L\'evy density $N$. 

\subsection{Bandwidth selection}\label{Bandwidth2}
 As $\beta$ is unknown, we need a data-driven selection of the bandwidth. We follow ideas given in  \cite{GL11} for density estimation.
We introduce a set of bandwidth of the form 
$H=  \lbrace  \frac{j}{M},  1 \leq j \leq M \rbrace$ with $M$ an integer to be specified later. Actually it is sufficient to control $\sum_{h\in H}h^{-w}$ for some $w$ so that more general set of bandwiths are possible.
We set: 
   \begin{eqnarray*}V(h)= C_0 \frac{\log(n\Delta)}{n h \Delta }
     \end{eqnarray*}
   with $C_0$ to be specified later. 
Note that $V(h)$ has the same order as the variance multiplied by $\log(n\Delta)$.
We also define  $\hat{g}_{h,h'}(x_0) = K_{h'} \star \hat{g}_h(x_0) = K_h \star \hat{g}_{h'}(x_0) $. This auxiliary estimator can also be written \begin{equation*}
\hat{g}_{h,h'}(x_0)= \frac{1}{n\Delta}\sum_{k=1}^n Z_k^{\Delta} K_{h'}\star  K_h(x_0-Z_k^{\Delta}).
\end{equation*} 
Lastly we set, as an estimator of the bias, 
\begin{equation*}\label{Ah}
A(h,x_0)= \sup_{h'\in H} \left[ |\hat{g}_{h,h'}(x_0) - \hat{g}_{h'}(x_0) |^2 - V(h')\right]_+.
\end{equation*} 
The adaptive bandwidth $h$ is chosen as follows:
\begin{eqnarray*}
\hat{h}=\hat{h}{(x_0)}\in\arg\min_{h \in H} \lbrace A(h,x_0) + V(h) \rbrace.
\end{eqnarray*}
We can state the following oracle inequality.

\begin{thm} \label{theomet2} 
We use a kernel satisfying $K1$ and  a set of bandwidth
$H=  \lbrace  \frac{j}{M},  1 \leq j \leq M \rbrace$ with $M= O((n\Delta)^{1/3})$.
Assume that $g$ satisfies G1, G2, G3(5)  and 
take \begin{equation}
      \label{cprime}
C_0=C_0(c)=\frac{c}{2\pi} \|K\|^2\left(\|(g^*)'\|_1+\|g^*\|_2^2\right)
     \end{equation}
with $c\geq 16\max(1,\|K\|_\infty).$
Then, for $\Delta\leq 1$,
\begin{eqnarray*}
\mathbb{E}[|g(x_0)- \hat{g}_{\hat{h}}(x_0)|^2] \leq C\left\lbrace \inf_{h \in H}\left\lbrace \|g- \mathbb{E}[\hat{g}_{h}] \|_\infty^2 + V(h)\right\rbrace  + \frac{\log(n\Delta)}{n\Delta}\right\rbrace
\end{eqnarray*}
\end{thm}

Thus our estimator $\hat{g}_{\hat{h}}$ has a risk as good as any of the collection $(\hat g_h)_{h\in H}$, up to a logarithmic term.

Note that the theorem is valid for $c$ large enough, say $c\geq c_0$. In the proof, we obtain the upper bound $16\max(1,\|K\|_\infty)$ for $c_0$, 
 unfortunately we can conjecture that this bound is not the optimal one.
 To obtain a sharper bound we have tuned $c_0$ in the simulation study.
 
The definition of the estimator uses $\|(g^*)'\|_1$ and $\|g^*\|_2^2$, but these quantities  can be estimated with a preliminar estimator of $g^*$.
%
%
%
%
More precisely, 
we set $K_{0}^{*}=\1_{[-1,1]}$ and
 $$\widehat{\|(g^*)'\|_1}=\int \left|\frac{1}{n\Delta}\sum_{k=1}^{n} (Z_k^{\Delta})^{2} K_{0}^{*}(uh_{1})e^{iuZ_k^{\Delta}}\right|du\quad\text{with }h_{1}=(n\Delta)^{-1/3},$$
 $$\widehat{\|g^*\|_2^2}=\|\hat g_{h_{2}}^*\|_2^2=\int \left|\frac{1}{n\Delta}\sum_{k=1}^{n} Z_k^{\Delta} K_{0}^{*}(uh_{2})e^{iuZ_k^{\Delta}}\right|^{2}du\quad\text{with }h_{2}=(n\Delta)^{-1/3}.$$
 We introduce the following regularity condition: a fonction $\psi$ belongs to the Sobolev space $Sob(\alpha)$ if
 $\int |\psi^{*}(u)|^{2}|u|^{2\alpha}du<\infty$.
Then, reinforcing the conditions on $g$,  we obtain a similar theorem with an empirical $C_0$. 
 \begin{thm} \label{theoEst} 
We use a kernel satisfying $K1$ and  $K2(\alpha)$ with $\alpha\geq1$, 
and  $M= O((n\Delta)^{1/3})$.
Assume that $g$ satisfies G1, G2, G3(32), G4(1), G5. Assume also that $g$ and $xg(x)$ belong to $Sob(1)$.
Take \begin{equation*}
     C_0=\frac{c}{2\pi} \|K\|^2\left(\widehat{\|(g^*)'\|_1}+\widehat{\|g^*\|_2^2}\right)
     \end{equation*}
with $c\geq 32\max(1,\|K\|_\infty).$
Then, for $n^{-1}\leq\Delta\leq C n^{-1/3}$,
\begin{eqnarray*}
\mathbb{E}[|g(x_0)- \hat{g}_{\hat{h}}(x_0)|^2] \leq C\left\lbrace \inf_{h \in H}\left\lbrace \|g- \mathbb{E}[\hat{g}_{h}] \|_\infty^2 + \E(V(h))\right\rbrace  + \frac{\log(n\Delta)}{n\Delta}\right\rbrace
\end{eqnarray*}
\end{thm}

Let us now conclude with the consequence of this theorem in term of rate of convergence. As already explained, as we need assumption G5 to control the bias, we can assume $\beta\geq 1$. Then $h_{opt}\propto (\log(n\Delta)/n\Delta)^{1/(2\beta+1)}\geq (n\Delta)^{-1/3}$ belongs to $H$ as soon as $M$ is larger than a constant times $(n\Delta)^{1/3}$. Hence we can state the following corollary. 

\begin{cor} 
 Assume that $g$ satisfies G1, G2, G3(5), G4($\beta$) with $\beta\geq 1$ and G5. We choose a kernel satisfying K1 and K2($\alpha$) with $\alpha\geq\beta$, and $M=\lfloor (n\Delta)^{1/3}\rfloor$. Take $C_0$ as in Theorem \ref{theomet2}
 (or as in Theorem~\ref{theoEst} with assumptions of this latter theorem).
Then, if $n^{-1}\ll\Delta\leq C n^{-1/3}$, 
\begin{eqnarray*}
\mathbb{E}[|g(x_0)- \hat{g}_{\hat{h}}(x_0)|^2] =O\left((\log(n\Delta)/n\Delta)^{-\frac{2\beta}{2\beta+1}}\right).
\end{eqnarray*}
\end{cor}

Then the price to pay to adaptivity is a logarithmic loss in the rate. Nevertheless this phenomenon is known to be unavoidable in pointwise estimation
(see \cite{BUC}). Thus $\hat{g}_{\hat{h}}(x_0)$ (resp. $\hat{g}_{\hat{h}}(x_0)/x_0$) is an adaptive estimator for $g(x_0)$ (resp. $N(x_0)$).

\section{Examples and Simulations}\label{Examples}

\begin{figure}[h!]
\begin{tabular}{ccc}
Ex 1 ($n\Delta=1000$) MISE$=0.032$ 
& Ex 1 ($n\Delta=2500$) MISE$=0.014$ \\ 
\includegraphics[scale=0.4]{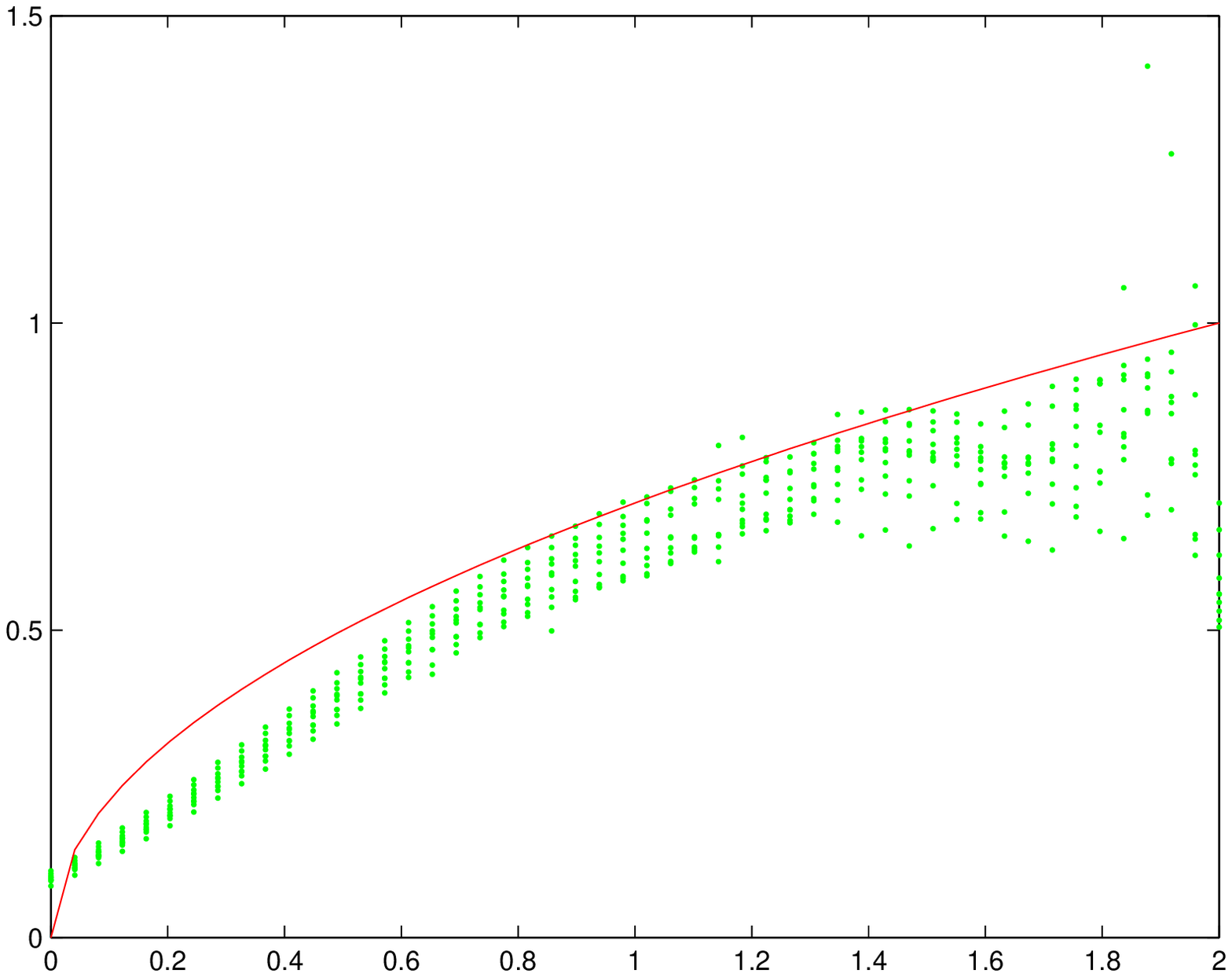}&\includegraphics[scale=0.4]{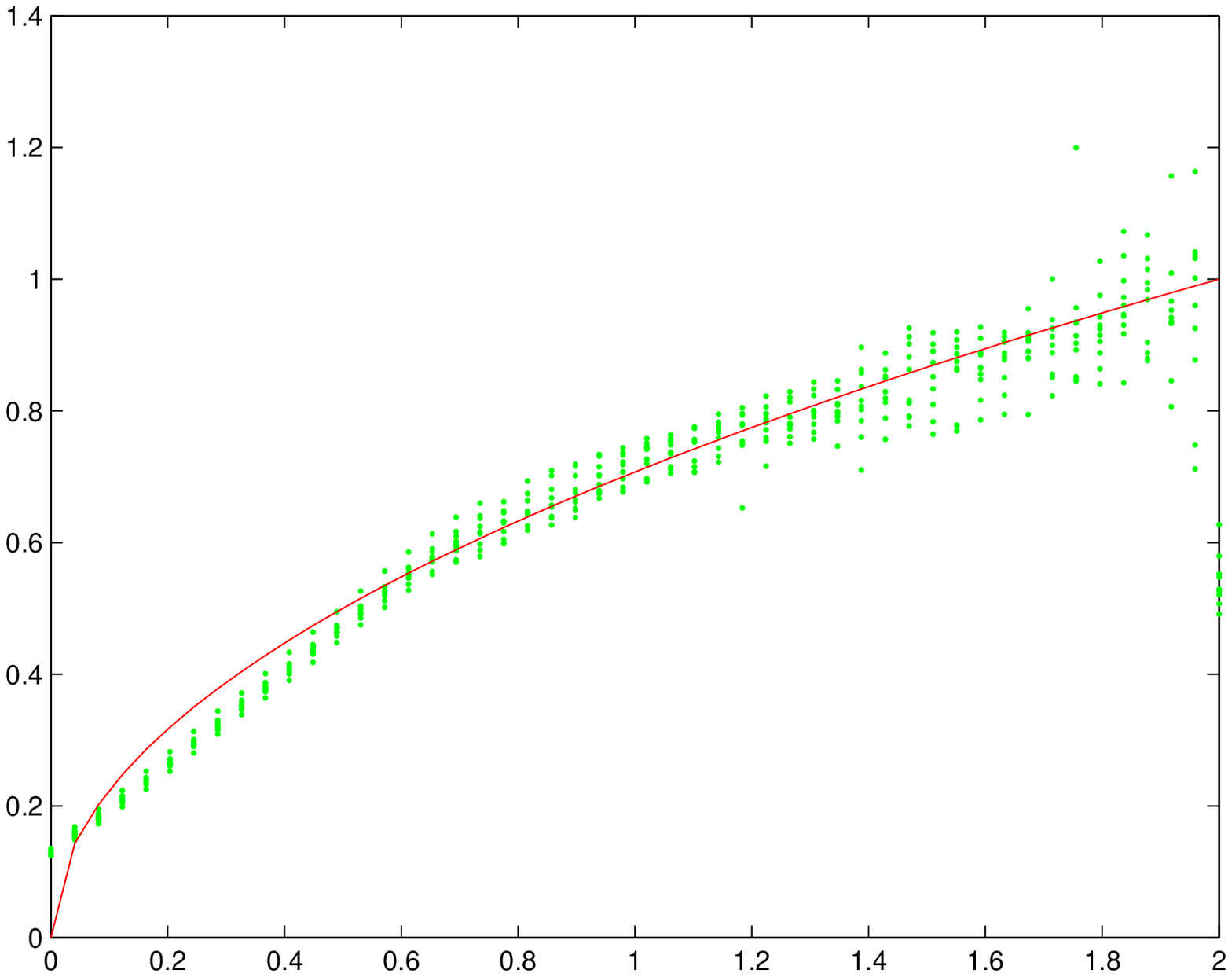}\\
&&\\
Ex 2 ($n\Delta=1000$) MISE$= 0.894$ 
& Ex 2 ($n\Delta=2500$) MISE$=0.057$ \\ 
\includegraphics[scale=0.4]{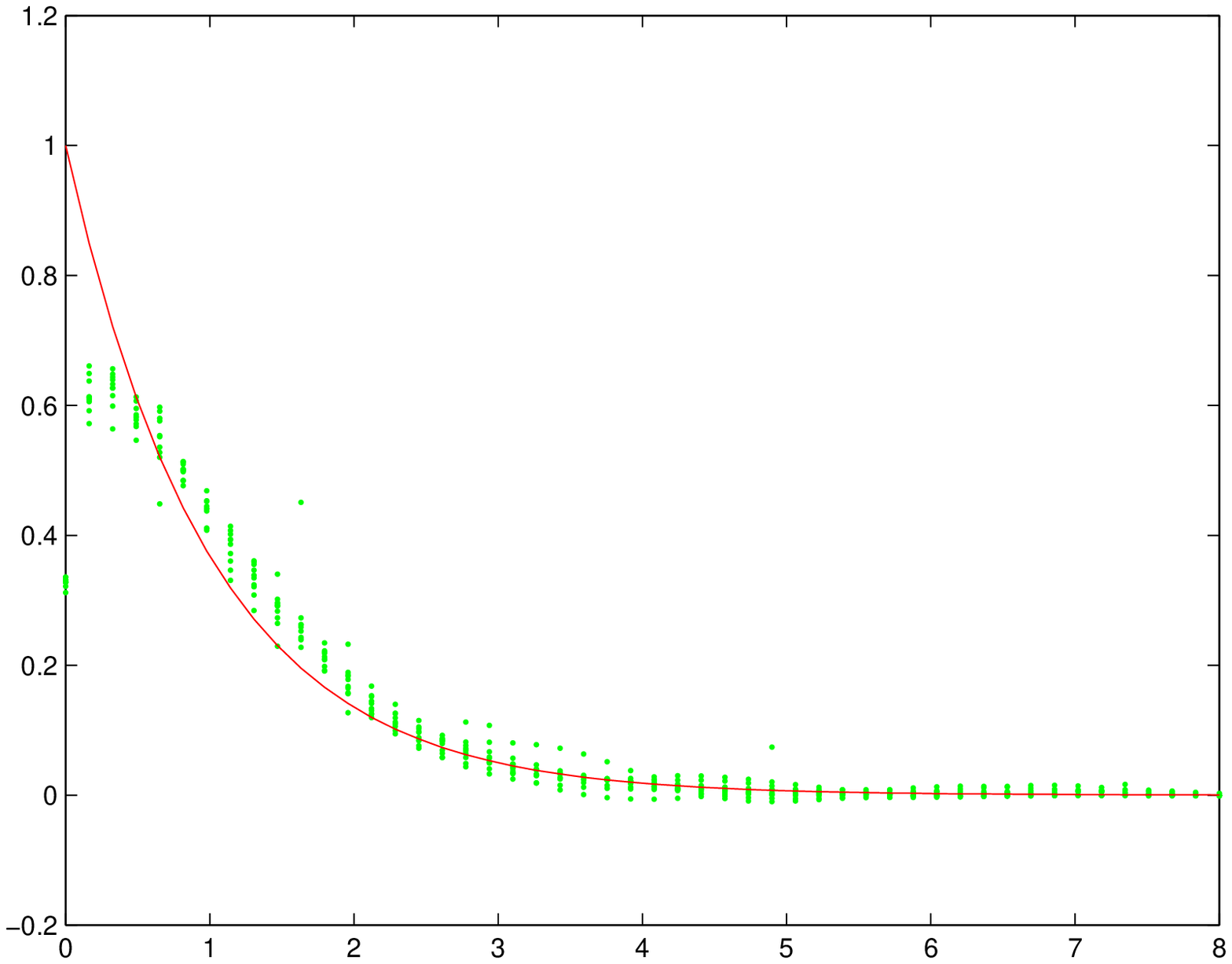}&\includegraphics[scale=0.4]{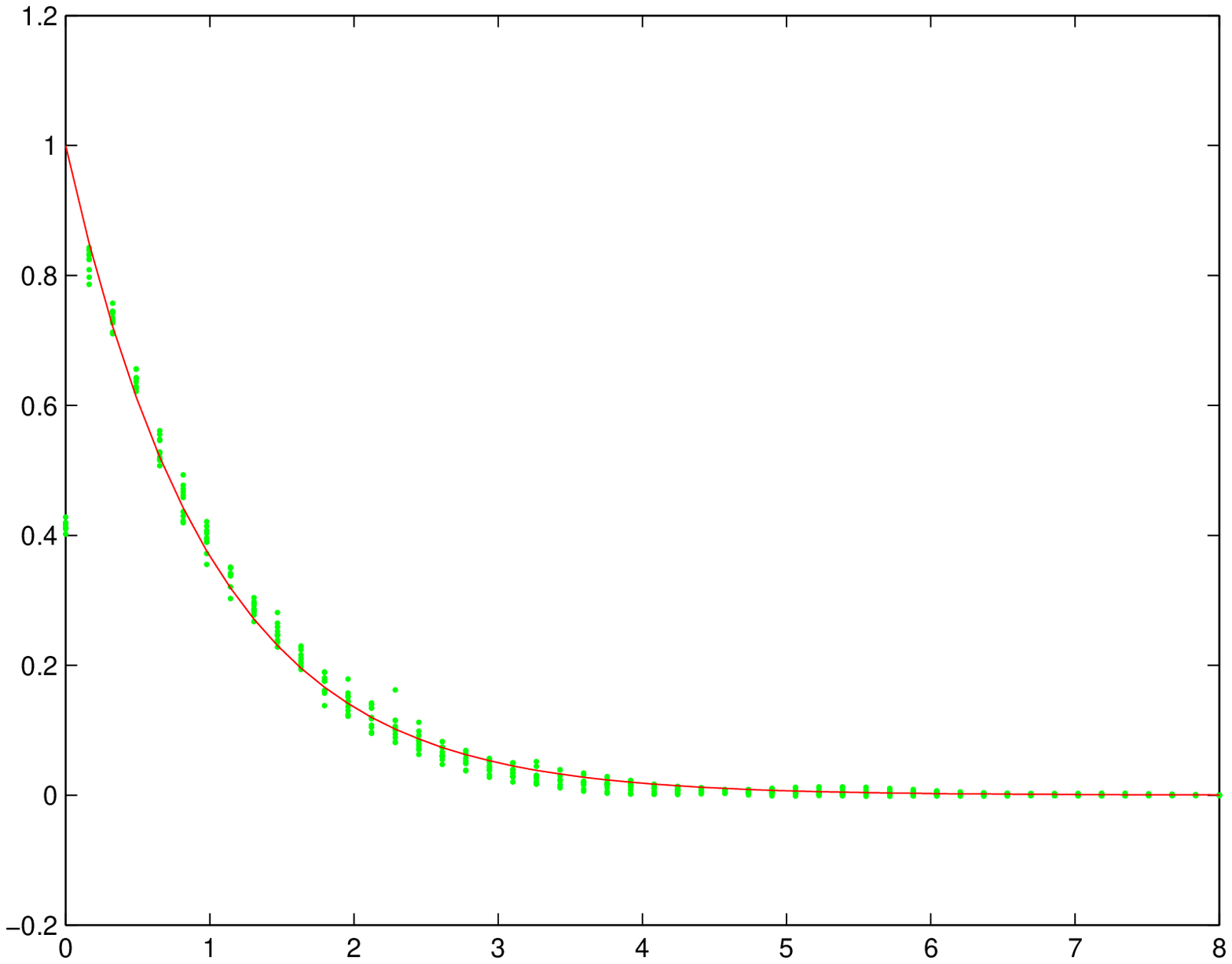}
\end{tabular}
\caption{Function $g$ (solid line) and estimators $\hat g_{\hat h}$ (dotted lines). }
\label{figure1}
\end{figure}

\begin{figure}[h!]
\begin{tabular}{ccc}
Ex 3 ($n\Delta=1000$) MISE$=0.009$ 
& Ex 3 ($n\Delta=2500$) MISE$=0.002$ \\
 \includegraphics[scale=0.4]{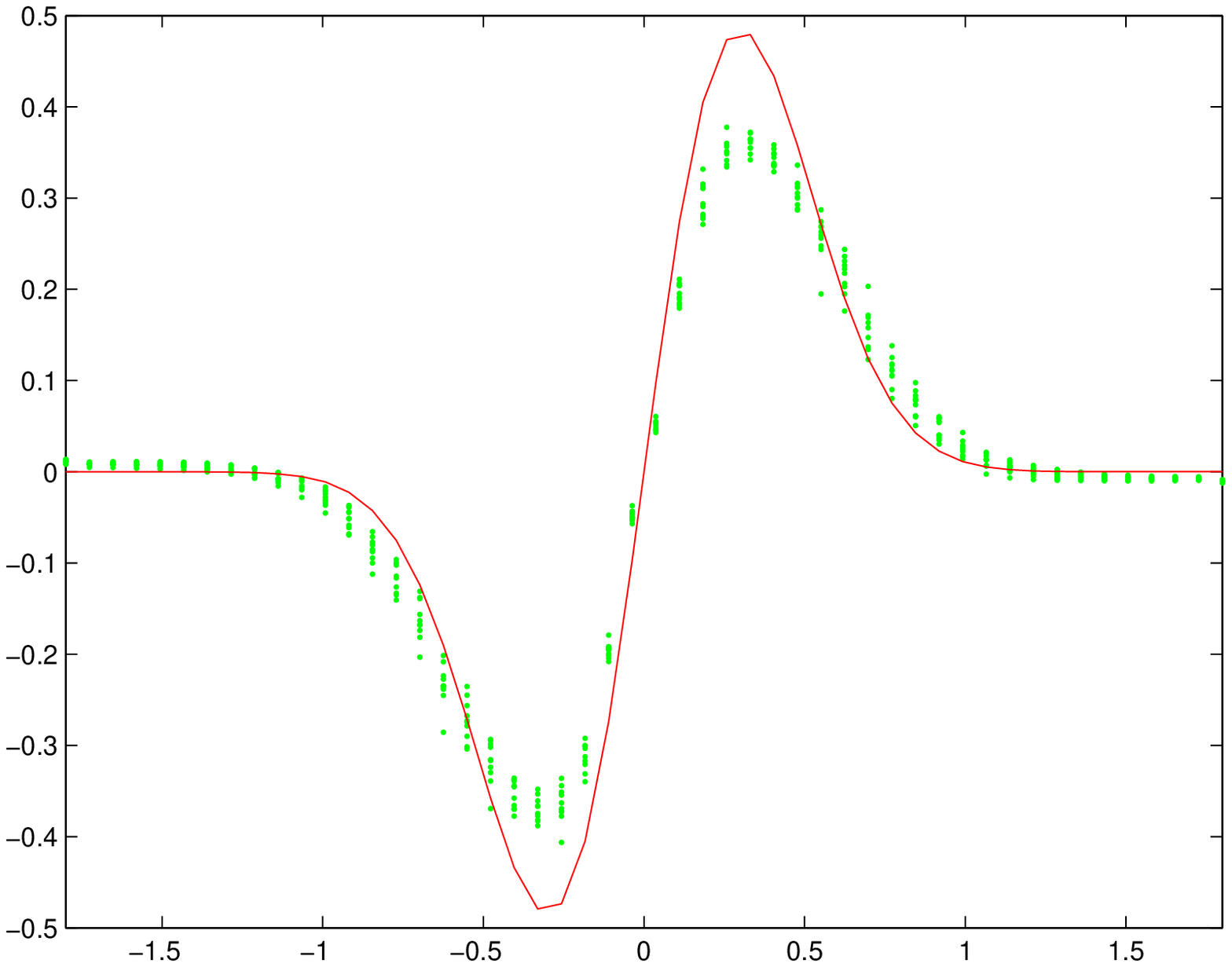}& \includegraphics[scale=0.4]{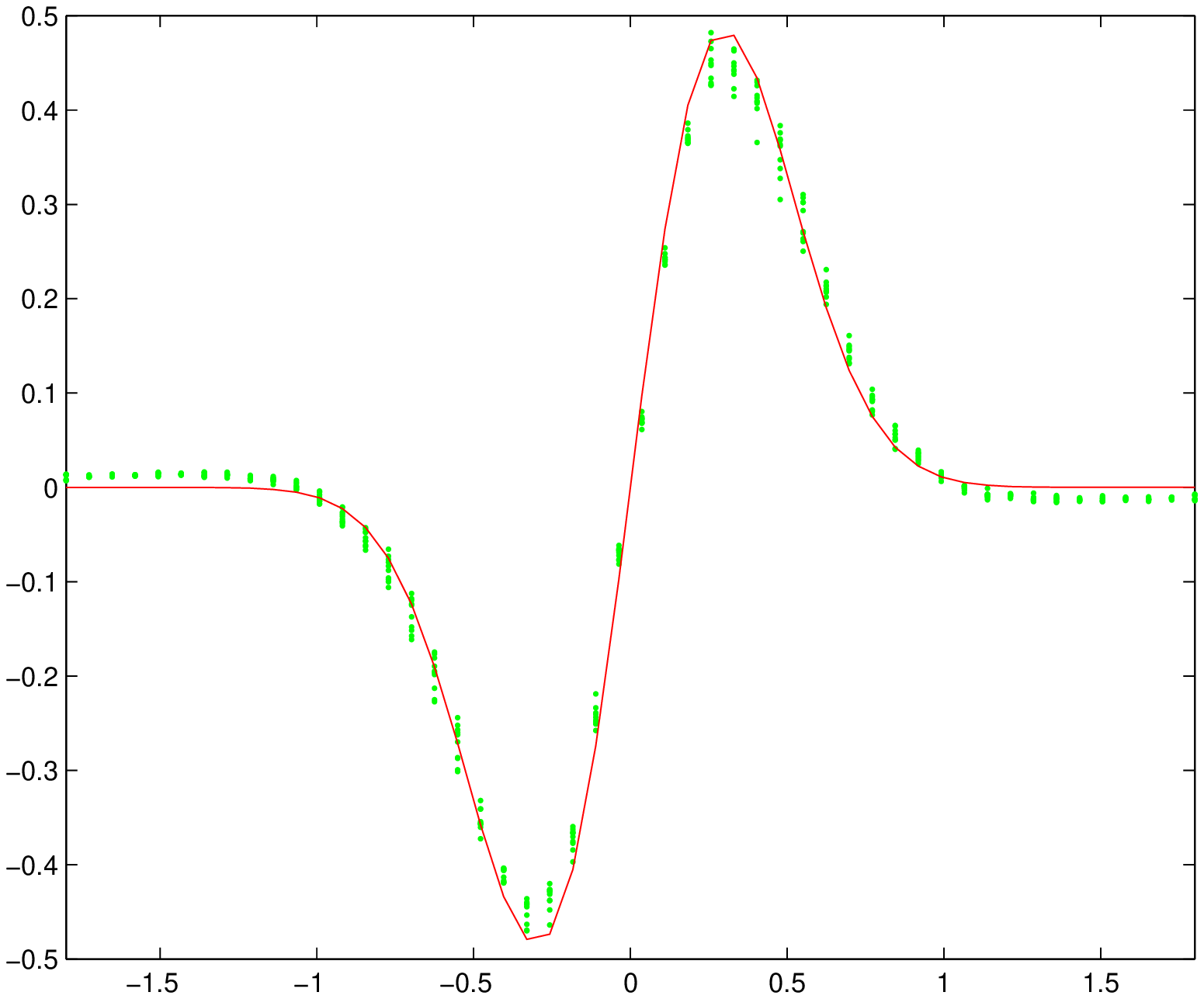}\\
&&\\
Ex 4 ($n\Delta=1000$) MISE$=0.811$ 
& Ex 4 ($n\Delta=2500$) MISE$=0.375$ \\
\includegraphics[scale=0.4]{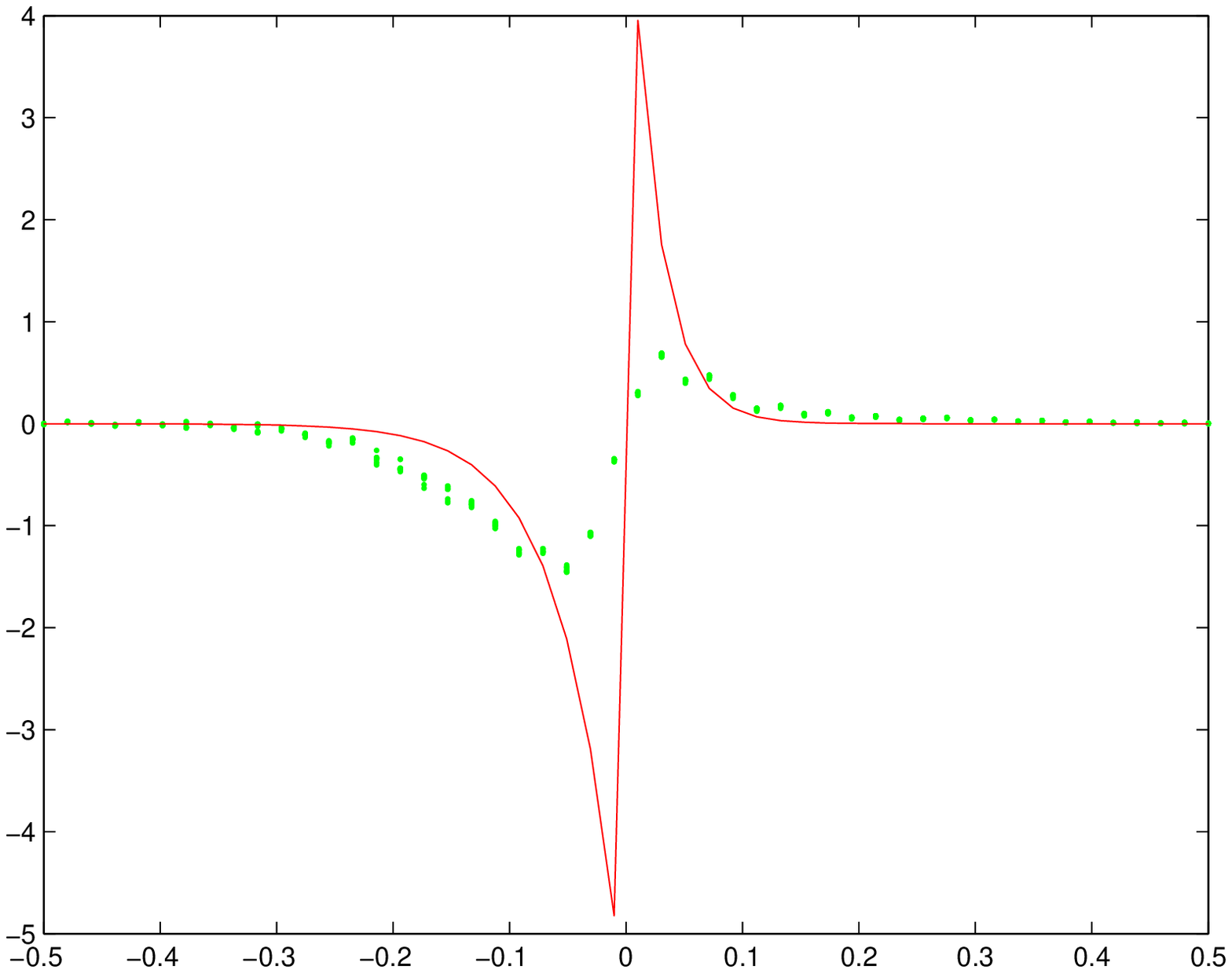}&\includegraphics[scale=0.4]{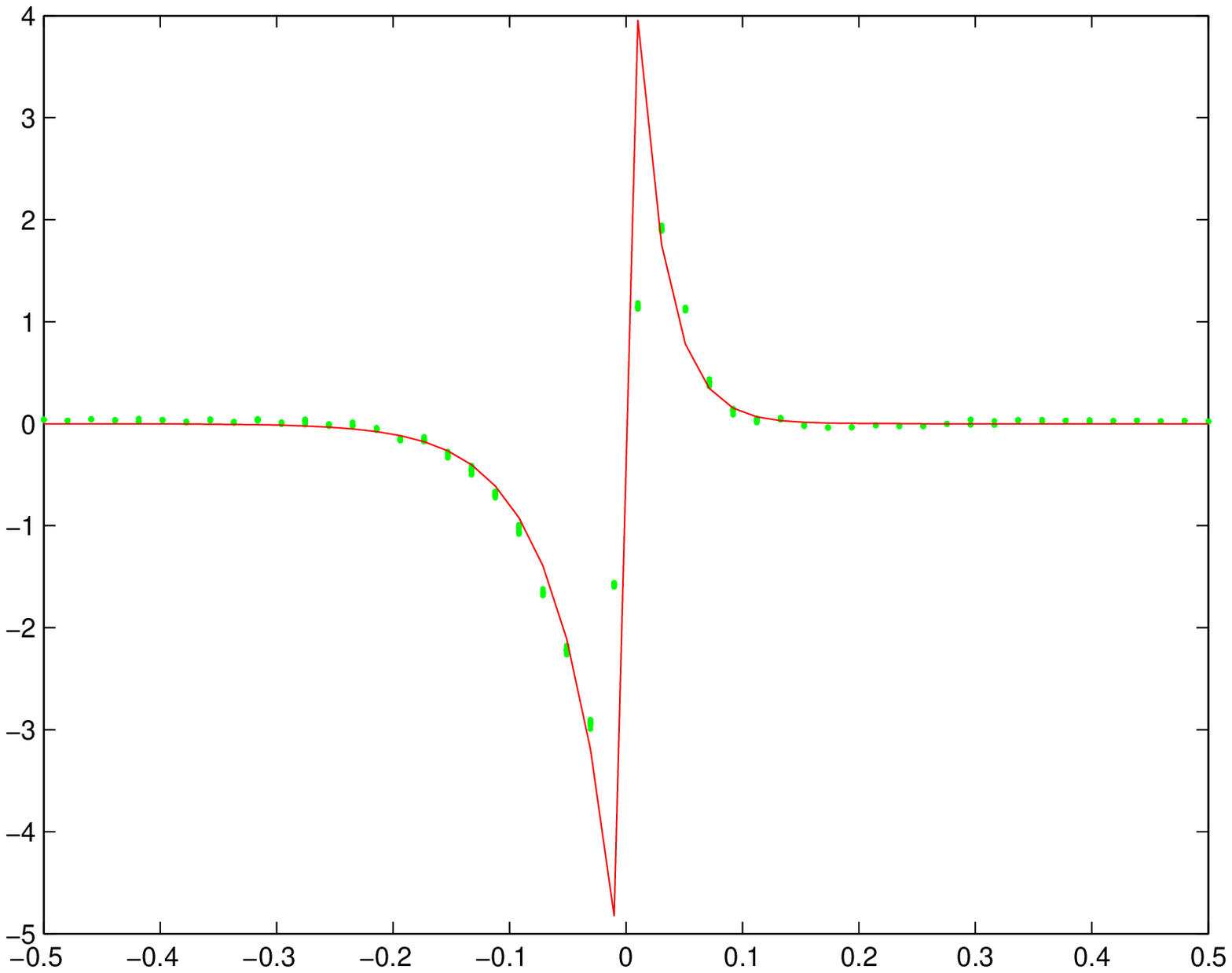}\\
\end{tabular}
\caption{Function $g$ (solid line) and estimators $\hat g_{\hat h}$ (dotted lines). }
\label{figure2}
\end{figure}

We have implemented the estimation method for four different processes (listed in Examples 1-4 below) with 
the kernel described in \eqref{noyauOP} (with $l=2$ and $u$ the Gaussian density).
The bandwidth set has been fixed to  $H=  \lbrace  \frac{j}{2M},  1 \leq j \leq M \rbrace$  with $M=\lfloor 2(n\Delta)^{-1/3}\rfloor$.
For the implementation, a difficulty is the proper calibration of the constant $c$ in (\ref{cprime}). This is usually done by a large number of preliminary simulations. We have chosen $c=0.1$ as the adequate value for a variety of models and number of observations. The estimation and adaptation are done for 50 points $x_{0}$ on the abscissa interval.
 For clarity, we have computed the Mean Integrated Square Error (MISE) of the estimators. Figures \ref{figure1}  and \ref{figure2} plot ten estimated curves corresponding to our four examples 
 with in the first column $\Delta=0.02, n=5.10^{3}$, and in the second $\Delta=0.05, n=5.10^{4}$.
This values of parameters can be interpreted as around hourly observations during few years.\\

\textbf{Example 1.} 
Let $L_t=\sum_{i=1}^{N_t}Y_i, $ where $(N_t)$ is a Poisson process with constant intensity $\lambda$ and $(Y_i)$ is a sequence of i.i.d random variables with density $f$ independent of the process $(N_t)$. Then, $(L_t)$ is a L\'evy process with characteristic function
\begin{equation}
\psi_{t}(u)= \exp\left( \lambda t \int_{\R}(e^{iux} -1) f(x) dx\right).
\end{equation}
Its L\'evy density is $N(x)= \lambda f(x)$ and thus $g(x)=\lambda xf(x)$. For our first example, we choose $\lambda=2$ and $f$ such that   
$g(x)=xf(x)=(1/2)\sqrt{x/2}$ for $0< x\leq 2$. Then assumption G4(1/2) holds (on $(0,2)$), but not G4($\beta$) for other $\beta$. Since $\beta$ is small, the rate of convergence is slow. The discontinuity in 2 damages the estimation as it can be seen in Figure~\ref{figure1}.

\textbf{Example 2.} 
Let $\alpha >0 $, $\gamma >0 $. The L\'evy-Gamma process $(L_t)$ with parameters $(\gamma,\alpha)$ is such that, for all $t>0$, $L_t $ has Gamma distribution with parameters $(\gamma t, \alpha)$, i.e the density:
\begin{eqnarray*}
\frac{\alpha^{\gamma t}}{\Gamma(\gamma t)} x^{\gamma t -1}e^{-\alpha x}\1_{x\geq 0}.
\end{eqnarray*}
The L\'evy density is $N(x)=\gamma x^{-1}e^{-\alpha x}\1_{x>0}$ so that $g(x)=\gamma e^{-\alpha x}\1_{x>0} $ satisfies assumptions G1, G2 and G3($p$). 
Here we choose $\alpha=\gamma=1$.
This example allows to study the role of the discontinuity in $0$, which invalidates assumptions G4-G5. We can observe that the estimation become very good if we move away from $0$.

\textbf{Example 3.} 
For our third example, we also choose a compound Poisson process, but with  $f$ the Gaussian density with variance $\delta^{2}.$ 
Thus $g(x)=\lambda xf(x)= \lambda xe^{-x^2/(2\delta^{2})}/(\delta\sqrt{2\pi}) $ and $g^*(u)= i\lambda\delta ue^{-\delta^{2}u^2/2}$.
Assumptions G1, G2, G3($p$),G5 hold for $g$. Moreover $g$ belongs to a H\"older class of regularity $\beta$ for all $\beta>0$. Thus the rate is close to  $(n\Delta/\log(n\Delta))^{-1}$, and the good performance of our estimator is visible on Figure~\ref{figure2}. 
Note that is the so-called Merton model used for describing the log price in financial modeling. Here we choose  $\lambda=2$ and $\delta=0.3$.

\textbf{Example 4.} 
Our last example is the Variance Gamma process, as described in \cite{MadanCarrChang98}.  
It is used for modeling the dynamics of the logarithm of stock prices. The process is obtained in evaluating a Brownian
motion at a time given by a L\'evy-Gamma process. Denoting $(B_{t})$ a standard Brownian motion, and $(X_{t})$ a  L\'evy-Gamma process with parameters $(1/\nu,1,\nu)$ independent of $(B_{t})$,
we set $L_{t}=\theta X_{t}+\sigma B_{X_{t}}$. Then $L_{t}$ is a L\'evy process, with 
$$g(x)=\frac{x\exp(\theta x /\sigma^{2})}{\nu|x|}\exp\left(-\frac1\sigma\sqrt{\frac2\nu+\frac{\theta^{2}}{\sigma^{2}}}|x|\right).$$
As in example 3, there is a discontinuity in $0$. 
Here we choose  $\theta=-0.1436$, $\sigma=0.1213$, $\nu=0.1686$: these are estimates of parameters for the S\&P index option prices studied in \cite{MadanCarrChang98}.


 \section{Irregular sampling}\label{IS}
 
For high frequency data, it is frequent that the sampling is irregular, i.e. the  interval $\Delta$ is not necessarily the same at each time. In this section we consider the following framework. The observations are $(L_{t_{k}},{k=1,...,n})$ where $(L_{t})$ is still a  L\'evy process with characteristic function
\eqref{fc}. For each $k\geq1$, we denote  $\Delta_{k}=t_{k}-t_{k-1}$ the sampling intervals. Notice that it includes the previous case when for each $k$, $\Delta_{k}=\Delta$.
 The increments are denoted by
 $Z_k= L_{t_k}- L_{t_{k-1}}$. In this context of irregular sampling, they are still independent but with non-identical distribution: $Z_{k}$ has the same law than
 $L_{\Delta_{k}}$. 
 To define an estimator, we observe that $\mathbb{E}\left[ Z_k e^{iuZ_k} \right] =\Delta_{k}\psi_{\Delta_{k}}(u) g^*(u)$,
 and then 
$$\mathbb{E}\left[\frac{1}{\sum_{k=1}^{n}\Delta_{k}}\sum_{k=1}^{n} Z_k e^{iuZ_k} \right] =
\left(\frac{\sum_{k=1}^{n}\Delta_{k}\psi_{\Delta_{k}}(u)}{\sum_{k=1}^{n}\Delta_{k}}\right) g^*(u).$$
Thus, denoting $\bar\Delta=\frac1n\sum_{k=1}^{n}\Delta_{k}$, we introduce 
  \begin{equation}\label{estimIS}
\hat g_h^*(u)=\frac{1}{n\bar\Delta}\sum_{k=1}^{n} Z_k e^{iuZ_k}K^*(hu) ,\quad
\hat g_h (x)= \frac{1}{n\bar\Delta} \sum_{k=1}^{n} Z_k K_h( x-Z_k)
\end{equation}
Additionally, for all real $\delta$, we denote $\overline{\Delta^{\delta}}=\frac1n\sum_{k=1}^{n}\Delta_{k}^{\delta}$. 
We can bound the Mean Squared Error of this estimate:

\begin{prop} \label{proprisque1IS}
Under G1, G2, G3(1), G4($\beta$), G5 and if $K$ satifies K1 and K2($\alpha$) with $\alpha\geq \beta$, we have
\begin{eqnarray} \label{borneIS}
MSE(x_0,h) \leq c_1 h^{2\beta} + c_2 \frac{1}{nh \bar\Delta} + {c'}_2 \frac{\overline{\Delta^2}}{nh\bar\Delta^{2}} + {c'}_1 \left(\frac{\overline{\Delta^2}}{\bar\Delta}\right)^{2}
\end{eqnarray}
with $c_1=  2 \left( { L}/{\lfloor \beta \rfloor !} \int |K(v)| |v|^\beta dv \right) ^2$, ${c'}_1= 2 (2\|g'\|_\infty{\| g \|}_1 {\| K \|}_1)^2$,
$c_2 = \|(g^*)'\|_1{\| K \|}^2_2 /(2\pi) $, ${c'}_2= {\| K \|}^2_2\| g \|^2_2 $.
\end{prop} 
The proof is similar to the case of regular sampling, therefore it is omitted. \\

In this section, we are still interested in the high frequency context: the asymptotic framework is $\bar\Delta\rightarrow 0$ and $n\bar\Delta\rightarrow \infty$ when $n\rightarrow \infty$.
We shall also assume that \begin{eqnarray} \label{condIS}
\frac{(\overline{\Delta^2})^2}{\bar\Delta}=O( n^{-1}) .
 \end{eqnarray}
 Condition (\ref{condIS}) is verified for instance if $\Delta_{k}=Ck^{-\alpha}$ with $\alpha\in[1/3,1]$. 
 Then we find the same rate of convergence replacing $\Delta$ by $\bar\Delta$:

\begin{prop} \label{proprisque2IS}
Under the assumptions of Proposition \ref{proprisque1IS} and under condition (\ref{condIS}), the choice $h_{opt}\propto((n\bar\Delta)^{-\frac{1}{2\beta + 1}})$  minimizes the risk bound (\ref{borneIS}) and gives $MSE(x_0,h_{opt})=O((n\bar\Delta)^{-\frac{2 \beta}{2\beta + 1}})$.
\end{prop}

As already noticed in \cite{CGCNeerlandica}, other estimation strategies than \eqref{estimIS} are possible.
 For each real  $\delta$, we obtain an estimator by setting 
 \begin{equation*}
\hat g_h (x)= \frac{1}{n\overline{\Delta^{\delta+1}}} \sum_{k=1}^{n} \Delta_{k}^{\delta}Z_k K_h( x-Z_k).
\end{equation*}
Under suitable conditions, this estimate has a MSE bounded by a constant times $(n\overline{\Delta^{\delta+1}}^{2}/\overline{\Delta^{2\delta+1}})^{-\frac{2 \beta}{2\beta + 1}}$. But, for all 
$\delta$, by the Schwarz inequality, $\overline{\Delta^{\delta+1}}^{2}/\overline{\Delta^{2\delta+1}}\leq \bar\Delta$. That is why we prefer estimator \eqref{estimIS}.

To build an adaptive estimator, we use the same method of bandwidth selection.
The set of bandwidth is still
$H=  \lbrace  \frac{j}{M},  1 \leq j \leq M \rbrace$. We also define 
\begin{equation*}
\hat{g}_{h,h'}(x_0)= K_{h'} \star \hat{g}_h(x_0) = \frac{1}{n\bar\Delta}\sum_{k=1}^n Z_kK_{h'}\star  K_h(x_0-Z_k)
\end{equation*} 
and we set as previously 
$A(h,x_0)= \sup_{h'\in H} \left[ |\hat{g}_{h,h'}(x_0) - \hat{g}_{h'}(x_0) |^2 - V(h')\right]_+$
with
   \begin{eqnarray*}V(h)= C_0 \frac{\log(n\bar\Delta)}{n h \bar\Delta }.
     \end{eqnarray*}
Then the estimator is $\hat{g}_{\hat{h}}(x_0)$ with 
$\hat{h}=\hat{h}{(x_0)}\in\arg\min_{h \in H} \lbrace A(h,x_0) + V(h) \rbrace.$

We can state the following oracle inequality (the proof is very similar to the one of Theorem~\ref{theomet2} and is therefore omitted). 
%

\begin{thm} \label{theoIS2} 
We use a kernel satisfying $K1$ and  $M= O((n\bar\Delta)^{1/3})$.
Assume that $g$ satisfies G1, G2, G3(5)  and 
take \begin{equation}
      \label{C0IS}
C_0=\frac{c}{2\pi} \|K\|^2\left(\|(g^*)'\|_1+\|g^*\|_2^2\right)
     \end{equation}
with $c\geq 16\max(1,\|K\|_\infty).$
Then, if $(\overline{\Delta^2})^2/\bar\Delta\leq 1$,
\begin{eqnarray*}
\mathbb{E}[|g(x_0)- \hat{g}_{\hat{h}}(x_0)|^2] \leq C\left\lbrace \inf_{h \in H}\left\lbrace \|g- \mathbb{E}[\hat{g}_{h}] \|_\infty^2 + V(h)\right\rbrace  + \frac{\log(n\bar\Delta)}{n\bar\Delta}\right\rbrace
\end{eqnarray*}
Moreover, if $g$ satisfies G5, G4($\beta$) with $\beta\geq 1$ and the kernel satisfying K1 and K2($\alpha$) with $\alpha\geq\beta$, and $M=\lfloor (n\bar\Delta)^{1/3}\rfloor$, $\bar\Delta\ll n^{-1}$ and $(\overline{\Delta^2})^2/\bar\Delta=O( n^{-1}) $, then
\begin{eqnarray*}
\mathbb{E}[|g(x_0)- \hat{g}_{\hat{h}}(x_0)|^2] =O\left((\log(n\bar\Delta)/n\bar\Delta)^{-\frac{2\beta}{2\beta+1}}\right).
\end{eqnarray*}

\end{thm}
Thus the rate of convergence in this case of irregular sampling is 
$(\log(n\bar\Delta)/n\bar\Delta)^{-\frac{2\beta}{2\beta+1}}$ provided that 
$(\overline{\Delta^2})^2/\bar\Delta=O( n^{-1})$.

\section{Proofs}\label{Proofs}

Let us first state two useful propositions (see Proposition 2.1 in \cite{CGC} and Proposition 2.1 in \cite{CG} for a proof).

\begin{prop} \label{proposable2}
Denote by $P_{\Delta}$ the distribution of $Z_{1}^{\Delta}$ and define $\mu_{\Delta}(dx)=\Delta^{-1}x P_{\Delta}(dx) $. 
If $\int_{\mathbb{R}} |x|N(x)<\infty$, the distribution $\mu_{\Delta} $ has a density $h_{\Delta} $ given by \begin{eqnarray*} h_{\Delta}(x)= \int g(x-y)P_{\Delta}(dy)= \mathbb{E}g(x-Z_1^{\Delta}).
\end{eqnarray*}
\end{prop}

\begin{prop} \label{proposable}
Let $p \ge 1$ an integer such that  $\int_{\R}|x|^{p-1} |g(x)| dx < \infty$. Then $\E(|Z_1^{\Delta}|^{p}) < \infty$ 
and $\E[(Z_{1}^{\Delta})^{p}]= \Delta\int_{\R} x^{p-1}g(x)dx  + o(\Delta)$. Moreover, if $g$ is integrable, $\mathbb{E}(|Z_{1}^{\Delta}|) \leq 2\Delta \| g \|_1$.
\end{prop}

\subsection{ Proof of Lemma \ref{monlem1}. }

First, we study $b_2(x_0)$ using Proposition~\ref{proposable2}:
\begin{eqnarray*}
b_2(x_0) &=& \frac{1}{h \Delta} \mathbb{E} \left[Z_1^{\Delta} K\left( \frac{x_0 - Z_1^{\Delta}}{h}\right)  \right] - \frac{1}{h} \int K\left(  \frac{x_0-u}{h}\right) g(u) du \\ 
&=& \frac{1}{h} \int K \left( \frac{x_0-u}{h} \right) \mathbb{E}[ g(u-Z_1^{\Delta}) - g(u) ] du.
\end{eqnarray*}
Now, applying the mean value theorem to $g$, we get
\begin{eqnarray*}
|b_2(x_0)| &=& \left| \frac{1}{h} \int K \left( \frac{x_0-u}{h} \right) \mathbb{E}[ -Z_1^{\Delta} g'(u_{Z_1}) ] du \right| 
\mbox{ with $u_{Z_1} \in [u-Z_1^{\Delta},u$ ] }\\ 
& \leq &  \| g' \|_{\infty} \| K \|_1 \mathbb{E}\left|Z_1^{\Delta}\right|  .
\end{eqnarray*}
From the results of Proposition \ref{proposable} we obtain 
\begin{eqnarray}\label{inegalite1}
|b_2(x_0)|  \leq & 2 \| g' \|_{\infty} \| K \|_1 \| g \|_1\Delta  .
\end{eqnarray}

To study $ b_1(x_0) = {K_h}\star g(x_0)-g(x_0) $, it is sufficient to use Taylor's theorem and $G4(\beta)$ (this is a classic computation, see \cite{TSY} for details) and we obtain 
\begin{eqnarray} \label{inegalite2} 
 | b_1(x_0)|  \leq   \frac{h^\beta L}{l!} \int |K(v)| |v|^\beta dv.
\end{eqnarray}

Gathering (\ref{inegalite1}) and (\ref{inegalite2})
completes the proof of Lemma \ref{monlem1}. 
$\Box$

\subsection{ Proof of Lemma \ref{monlem2}. }
As
the $Z_k^{\Delta}$ are i.i.d., we have:
\begin{eqnarray*}
\Var[\widehat{g}(x_0)] = \Var\left[ \frac{1}{nh\Delta} \sum_{k=1}^{n} Z_k^{\Delta} K\left(  \frac{x_0-Z_k^{\Delta}}{h}\right) \right]  = \frac{1}{n(h\Delta)^2} \Var\left[  Z_1^{\Delta} K\left(  \frac{x_0-Z_1^{\Delta}}{h}\right) \right].
\end{eqnarray*}
Thus,
\begin{eqnarray*}
\Var[\widehat{g}(x_0)]   \leq  \frac{1}{n(h\Delta)^2} \mathbb{E}\left[  (Z_1^{\Delta})^2 K^2\left(  \frac{x_0-Z_1^{\Delta}}{h}\right) \right].
\end{eqnarray*}
Writing
\begin{eqnarray*}
 K^2\left(  \frac{x_0-Z_1^{\Delta}}{h}\right)  = {\left| \frac{1}{2 \pi} \int K^*(u)e^{-i \frac{(x_0-Z_1^{\Delta})u}{h}} du \right| }^2,
\end{eqnarray*}
we obtain with $v=u/h $
\begin{eqnarray*}
\Var[\widehat{g}(x_0)]  & \leq  &  \frac{1}{n \Delta ^2 } \mathbb{E}\left[  (Z_1^{\Delta})^2 {\left| \frac{1}{2 \pi} \int {{K}^{*}}(vh)e^{-i (x_0-Z_1^{\Delta})v} dv \right|}^2 \right] \\ & \leq & \frac{1}{n \Delta ^2 {(2 \pi)}^2} \mathbb{E}\left[\iint  Z_1^{\Delta}e^{iZ_1^{\Delta}v} {{K}^{*}}(vh)e^{-ix_0v}  \overline{Z_1^{\Delta}e^{iZ_1^{\Delta}u} {{K}^{*}}(uh)e^{-ix_0u} }dvdu \right].
\end{eqnarray*}
Using Fubini and $\mathbb{E}[(Z_1^{\Delta})^2 e^{iZ_1^{\Delta}(v-u)}]= -\psi''_{\Delta}(v-u) $  we find
\begin{eqnarray*}
\Var[\widehat{g}(x_0)]  & \leq  & \frac{1}{n \Delta ^2 {(2 \pi)}^2} \iint |-  {\psi_\Delta}''(v-u)  {{K}^{*}}(vh) {{K}^{*}}(uh)| dvdu 
\end{eqnarray*}
Now the following formula
\begin{eqnarray*} \label{relation2}
{\psi_{\Delta}}''= i\Delta {\psi_{\Delta}}' g^* + i\Delta \psi_{\Delta} {g^*}' = -{\Delta}^2 {\psi_{\Delta}} {g^{*}}^2 + i \Delta \psi_{\Delta} {g^*}'.
\end{eqnarray*}
gives $ \Var[\widehat{g}(x_0)] \leq T_1 +T_2$ with 
\begin{eqnarray*}
T_1=\frac{1}{n \Delta ^2 {(2 \pi)}^2} \iint  | {\Delta}^2 {\psi_\Delta}(v-u)(g^*)^2 (v-u)  {{K}^{*}}(vh) {{K}^{*}}(uh)| dvdu \\
T_2=\frac{1}{n \Delta ^2 {(2 \pi)}^2} \iint  |\Delta {\psi_\Delta}(v-u)(g^*)' (v-u)  {{K}^{*}}(vh) {{K}^{*}}(uh)| dvdu.
\end{eqnarray*}
We first bound $T_2$:
\begin{eqnarray*}
T_2 
& \leq & \frac{1}{n \Delta {(2 \pi)}^2} \sqrt{\iint  |{\psi_\Delta}(v-u)| |(g^*)' (v-u)| {| {{K}^{*}}(vh)|}^2  dvdu} \\ && \hspace{2cm} \times  \sqrt{\iint |{\psi_\Delta}(v-u)||(g^*)' (v-u)| {|{{K}^{*}}(uh)|}^2 dvdu } \\ & \leq & \frac{1}{n \Delta {(2 \pi)}^2} \int {| {{K}^{*}}(vh)|}^2 dv \int |{\psi_\Delta}(z)||(g^*)' (z)|   dz\\ 
& \leq & \frac{1}{n h  \Delta {(2 \pi)}^2} \int {| K^*(u)|}^2 du \int |(g^*)' (z)|   dz , \mbox{ because $|{\psi_\Delta}(z)|\leq 1 $}\\& \leq & \frac{ \| K \|_2^2}{2 \pi n h  \Delta}\int |(g^*)' (z)| dz 
\end{eqnarray*}
 where $(g^*)'$ exists and is integrable by G2. 
Following the same line for the study of $T_1$, we get
\begin{eqnarray*}
T_1 \leq \frac{{ {{\| K \|}^2_2}}}{ 2 \pi n h }\int |(g^*)^2 (z)| dz \leq \frac{{{\| K \|}^2_2 } {{\| g \|}^2_2 }}{n h },
\end{eqnarray*} 
This completes the proof of Lemma \ref{monlem2}. $\Box$

\subsection{Proof of the lower bound}\label{preuveborneinf}
Here we prove 
Theorem~\ref{lower}
The essence of the proof is to build two functions $g_0$ and $g_1$ which are far in term of pointwise distance but with close associated distribution.
Let $$g_{0}(x)=xf_\lambda(x)=\frac1\pi\frac{\lambda x}{1+(\lambda x)^2}$$ 
where $f_\lambda$ is the density of the Cauchy distribution $C(0,\lambda)$ with scale parameter $\lambda$.
Here $\lambda$ is a positive and small enough real (it will be made precise later).
Now let $K$ a infinitely differentiable and even function such that $\int K=0$, $K(0)\neq 0$ and
$K(x)=|x|^{-2}$ for $|x|$ large enough (say for $|x|> B$). Using this auxiliary function $K$, we can define 
$$g_1(x)=g_{0}(x)+ ch_n^\beta K\left(\frac{x-x_0}{h_{n}}\right)x$$
where $c$ is a constant to be specified later and $$h_n=(n\Delta)^{-\frac{1}{2\beta+1}} .$$
We denote $N_0(x)=g_0(x)/x$ and $N_1(x)=g_1(x)/x$. 
Remark that if $L_{0,t}=\sum_{i=1}^{N_t}Y_i$ is a compound Poisson process with $N_t$ a Poisson process of intensity 1 and 
$Y_i$ Cauchy $C(0,\lambda)$ variables, then its characteristic function is 
$$\psi_{0,t}(u)=\exp{(t \int_{{\mathbb R}}(e^{iux} -1) N_0(x)dx)}$$ and 
$Z_k^{0,\Delta}=L_{0,k\Delta}-L_{0,(k-1)\Delta}$ has distribution $P_{0}(dx)=e^{-\Delta}\delta_0(dx)+\varphi_0(x)dx$
with $$\varphi_0(x)=\sum_{k=1}^\infty e^{-\Delta}\frac{\Delta^{k}}{k!} f_\lambda^{*k}(x).$$
Moreover $N_1$ is a density. Indeed 
the definition of $K$ guarantees that  $\int N_1(x)dx=
\int N_0(x)dx + ch_n^\beta \int K\left(\frac{x-x_0}{h_{n}}\right)dx=1$.
And to ensure the positivity of $N_1$, it is sufficient to prove that
$|N_1-N_0|\leq N_0$. But, if $|x|>|x_0|+Bh_n$ ,
$$N_0^{-1}(x)|N_1(x)-N_0(x)|\leq C c h_n^{\beta+2}  x^{2}|x-x_0|^{-2}\leq 1$$
for $c$ small enough, and if $|x|\leq |x_0|+ B h_n$, 
$$N_0^{-1}(x)|N_1(x)-N_0(x)|\leq C c h_n^{\beta} (1+(\lambda (|x_0|+Bh_n))^2)\|K\|_\infty\leq 1$$
for $c$ small enough. 
Then, if $L_{1,t}=\sum_{i=1}^{N_t}Y_i$ with $N_t$ a Poisson process of intensity 1 and 
$Y_i$ random variables with density $N_1$, it is a L\'evy process with L\'evy measure $N_1(x)dx$.
We denote $\psi_{1,\Delta}$ the characteristic function of $L_{1,\Delta}$ with distribution $P_1$, and 
$\varphi_1$ the function such that $P_{1}(dx)=e^{-\Delta}\delta_0(dx)+\varphi_1(x)dx$.

Now let us denote for two probability measures $P$ and $Q$, $\chi^2(P,Q)=\int \left(dP/dQ-1\right)^2 dQ$. 
In the sequel we show that
\begin{itemize}
\item[1)]  $g_0,g_1$ belong to ${\mathcal H}(\beta,L)$,
\item[2)] $|g_1(x_0)-g_0(x_0)|\geq C(n\Delta)^{-\frac{\beta}{2\beta+1}}$,
\item[3)]  $\chi^2(P_1^n,P_{0}^n)\leq C<\infty$ where $P_{{1}}^n$ (resp. $P_{0}^n$) is the distribution of a sample $Z_1^\Delta, \dots, Z_n^\Delta$ s.t the associated L\'evy process $L_{0}$ (resp. $L_1$) has L\'evy measure $N_0(x)dx$ (resp. $N_1(x)dx$).

\end{itemize}
Then it is sufficient to use Theorem 2.2 (see also p.80) in \citet{TSY} to obtain Theorem \ref{lower}.
In the following we denote all constants by $C$, even if it changes from line to line.\\

\noindent Proof of 1).  {\it Belonging to the H\"{o}lder space } \\
To prove that our hypotheses belong to $\mathcal{H}(\beta,L)$,  it is sufficient to show that, for $i=0,1$,
$\|g_{i}^{(k+1)}\|_p\leq L$ where $k=\lfloor \beta\rfloor$ and $p^{-1}=1+k-\beta$. Indeed H\"older inequality gives
$$|g_i^{(k)}(x)-g_i^{(k)}(y)|=
\left|\int g_{i}^{(k+1)}(v) \1_{[x,y]}(v)dv\right|
\leq \|g_{i}^{(k+1)}\|_p|x-y|^{\beta-k}\quad\text{ for all } x, y.$$ 
When $x$ goes to infinity, $g_0^{(k+1)}(x)=C\lambda^{-1} x^{-k-2}+o(x^{-k-2})$ so it belongs to 
$\mathbb{L}^p$ since $p(k+2)=(k+2)/(k+1-\beta)>1$. Choosing $\lambda$ small enough ensures 
$\|g_{0}^{(k+1)}\|_p\leq L/2\leq L$.

Now to study $g_1$, we can write
\begin{eqnarray*}
(g_1-g_0)^{(k+1)}(x)=c xK^{(k+1)}\left(\frac{x-x_0}{h_n}\right)h_n^{\beta-k-1}+c(k+1) K^{(k)}\left(\frac{x-x_0}{h_n}\right)h_n^{\beta-k}.
\end{eqnarray*}
Let us see if this two terms are in $\mathbb{L}^p.$ Writing $x=x-x_0+x_0$ and changing variables
\begin{eqnarray*}
\int  \left|xK^{(k+1)}\left(\frac{x-x_0}{h_n}\right)\right|^pdx \leq 2^{p-1} h_n^{p+1}\int |v K^{(k+1)}(v)|^p dv 
+2^{p-1}|x_0|^ph_n\int |K^{(k+1)}(v)|^p dv .
\end{eqnarray*}
These integrals are finite since $vK^{(k+1)}(v)=v^{-(2+k)}$ for $v$ large enough and $p(k+2)=(k+2)/(k+1-\beta)>1$. In the same way 
\begin{eqnarray*}
\int  \left|K^{(k)}\left(\frac{x-x_0}{h_n}\right)\right|^pdx \leq h_n\int |K^{(k)}(v)|^p dv .
\end{eqnarray*}
Thus \begin{eqnarray*}
\|(g_1-g_0)^{(k+1)}\|_p^p\leq C c^p  ( h_n h_n^{p(\beta-k-1)} +h_n h_n^{p(\beta-k)})
\leq C c^p h_n^{p(1/p+\beta-k-1)}\leq C c^p \leq (L/2)^p
\end{eqnarray*}
for suitable $c$. 
Then $g_1-g_0$ belongs to ${\mathcal H}(\beta, L/2)$ and $g_1$ belongs to ${\mathcal H}(\beta, L)$.\\

\noindent Proof of 2). {\it Rate} \\
By assumption, $x_0\neq 0$ and we can see that
$|g_1(x_0)-g_0(x_0)|= ch_n^\beta|K(0)x_0|$ 
with $K(0)\neq 0$. Since   $h_n=(n\Delta)^{-\frac{1}{2\beta+1}} $, this quantity  has the announced order of the rate: $(n\Delta)^{-\frac{\beta}{2\beta+1}}$ .\\

\noindent Proof of 3). {\it Chi-square divergence}\\
Since the observations are i.i.d., $\chi^2(P_1^n,P_0^n)=(1+\chi^2(P_1,P_0))^n-1$. 
Thus, it is sufficient to prove that
$\chi^2(P_1,P_0)=O(n^{-1})$  where
$$\chi^2(P_1,P_0)=\int_{x\neq 0} \left(\frac{\varphi_1(x)}{\varphi_0(x)}-1\right)^2  \varphi_0(x)dx.$$
Indeed $P_1(\{0\})=e^{-\Delta}=P_0(\{0\})$.
Now let us remark that for $n$ large enough 
$$\varphi_0(x)=\sum_{k=1}^\infty e^{-\Delta}\frac{\Delta^{k}}{k!} f_\lambda^{*k}(x)\geq 
e^{-\Delta}\Delta f_\lambda(x)\geq \Delta e^{-C}\lambda\pi^{-1}/(1+(\lambda x)^2)$$
since $\Delta$ is bounded. 
Then $\varphi_0(x)\geq C^{-1} \Delta x^{-2}$ for $|x|$ large enough, say $|x|\geq A $
and  $\varphi_0(x)\geq C^{-1}\Delta $ for $|x|\leq A $.
Next we write 
$\chi^2(P_1,P_0)= \int_{x\neq 0} \left(\varphi_1(x)-\varphi_0(x)\right)^2 (\varphi_0(x))^{-1}dx=I_1+I_2$ where $I_1$ is the integral for $|x|< A$ and $I_2$ for $|x|\geq A$. We will bound these two terms separately.

Since $\varphi_0(x) \geq C^{-1}\Delta $ for $|x|$ small 
\begin{eqnarray}\nonumber
I_1&=& \int_{|x|< A}\!\!\left( \varphi_1(x)-\varphi_0(x)\right)^2(\varphi_0(x))^{-1}dx
\leq  C\Delta^{-1}\!\!\int_{|x|< A}\!\!\left(\varphi_1(x)-\varphi_0(x)\right)^2dx. 
\end{eqnarray}
For $i=0,1$,
the Fourier tranform of $\varphi_i$ is 
 $\psi_{i,\Delta}(u)-P_i(\{0\})$. 
Thus Parseval equality gives 
\begin{eqnarray*}
 I_1\leq C \Delta^{-1}\int \left|\psi_{1,\Delta}(u)-\psi_{0,\Delta}(u)\right|^2du.
\end{eqnarray*}
In order to get a bound on $|\psi_{1,\Delta}-\psi_{0,\Delta}|$, we apply the mean value theorem:
\begin{eqnarray*}
 |\psi_1(u)-\psi_0(u)|\leq \sup_{z\in I_u} |e^z| |\Delta\int (e^{iux}-1)(N_1(x)-N_0(x))dx|
\end{eqnarray*}
where $I_u$ is the segment in $\mathbb{C}$ between $a_u=\Delta\int (e^{iux}-1)N_0(x)dx$ and $b_u=\Delta\int (e^{iux}-1)N_1(x)dx$.
But 
\begin{eqnarray*}
 \int (e^{iux}-1)(N_1(x)-N_0(x))dx=ch_n^\beta\int (e^{iux}-1)K\left(\frac{x-x_0}{h_n}\right)dx
=ch_n^{\beta+1}e^{iux_0}K^*(h_nu).
\end{eqnarray*}
Note that this quantity is well defined since $K$ belongs to $\mathbb{L}^1$. 
Thus 
\begin{eqnarray*}
 |\psi_1(u)-\psi_0(u)|\leq & ( \sup_{z\in I_u} e^{\mathfrak{R}(z)} )\Delta c h_n^{\beta+1}|K^*(h_nu)|
\end{eqnarray*}
where $\mathfrak{R}(x)$ means the real part of $x$.
We can compute $\mathfrak{R}(a_u)=a_u=\Delta (N_0^*(u)-1)=\Delta(\exp(-|u/\lambda|)-1)\leq 0$ and 
$$\mathfrak{R}(b_u)=\mathfrak{R}(\Delta (N_0^*(u)-1+(N_1-N_0)^*(u)))=\Delta (N_0^*(u)-1+ch_n^{\beta+1}\mathfrak{R}(K^*(h_nu)e^{iux_0})).$$
Since $K$ is even, 
\begin{eqnarray*}
\mathfrak{R}(b_u)=\Delta (\exp(-|u/\lambda|)-1+ch_n^{\beta+1}K^*(h_nu)\cos(ux_0))
\leq c\Delta h_n^{\beta+1}\|K^*\|_\infty\leq C 
\end{eqnarray*}
so that 
\begin{eqnarray}\label{bornepsi}
 |\psi_1(u)-\psi_0(u)|\leq& e^C\Delta c h_n^{\beta+1}|K^*(h_nu)|.
\end{eqnarray}
Then
\begin{eqnarray}
 I_1 \leq C \Delta^{-1}\int \left|\Delta h_n^{\beta+1}K^*(h_nu)\right|^2du\leq C \Delta h_n^{2\beta+1}.  \label{I1}
\end{eqnarray}

Let us now bound the term $I_2$, using that
$\varphi_0(x) \geq C^{-1}\Delta x^{-2} $ for $|x|$ large enough
\begin{eqnarray*}
I_2&= &\int_{|x|\geq A} \frac{\left(\varphi_1(x)-\varphi_0(x)\right)^2}{\varphi_0(x)}dx
 \leq C \Delta^{-1} \int \left(\varphi_1(x)-\varphi_0(x)\right)^2x^2dx.
 \end{eqnarray*}
But $F=\varphi_1-\varphi_0$ has Fourier transform 
$$F^*=\psi_{1,\Delta}-\psi_{0,\Delta}=\exp(\Delta(e^{-|u/\lambda|}+ ch_n^{\beta+1}K^*(h_nu)e^{iux_0}-1))
-\exp(\Delta(e^{-|u/\lambda|}-1))$$ 
and this function is differentiable everywhere exept at $u=0$, with derivative
$$F^{*\prime}=\Delta\gamma_1\psi_{1,\Delta}-\Delta\gamma_0\psi_{0,\Delta}$$
where 
$$\gamma_0(u)= -{\rm  sign}(u). e^{-|u/\lambda|}/\lambda,\qquad \gamma_1(u)=\gamma_0(u)+ch_n^{\beta +1}e^{iux_0}(ix_0K^*(h_n u)+h_nK^{*\prime}(h_n u)).$$
Let us now prove that the Fourier transform of $F^{*\prime}$ is $-2\pi ixF(-x)$.
Let us write  the factorization
\begin{equation}\label{decompo}
\Delta^{-1}F^{*\prime}=\gamma_1\psi_{1,\Delta}-\gamma_0\psi_{0,\Delta}=(\gamma_1-\gamma_0)\psi_{1,\Delta}+\gamma_0(\psi_{1,\Delta}-\psi_{0,\Delta})
\end{equation} with  $|\psi_{1,\Delta}|\leq 1$.
Since $K^*$ and $K^{*\prime}$ are uniformly bounded,  $\gamma_1-\gamma_0$ is bounded as well. In the same way, the inequality \eqref{bornepsi} entails that $\|\psi_{1,\Delta}-\psi_{0,\Delta}\|_\infty<\infty$, so that $F^{*\prime}$ is bounded. Thus 
$F^*$ is Lipschitz and absolutely continuous. Moreover, using again \eqref{decompo}, we can see that $F^{*\prime}$ is integrable
(we can choose $K$ such that $K^*$ is integrable, for example take for $K$ the difference between the Cauchy density and the normal density).
Then, according to \cite{Rudin}, the Fourier transform of $F^{*\prime}$ is $-ixF^{**}(x)$ (it is in fact a simple integration by parts). 
Since $F^*$ is integrable,  $F^{**}(x)=2\pi F(-x)$ almost everywhere, and we have proved that  $(F^{*\prime})^*(x)=-2\pi ixF(-x)$ a.e..
%
Next, 
the Parseval equality provides $\int |xF(x)|^2dx=(2\pi)^{-1}\int |F^{*'}(u)|^2du$. Thus 
\begin{eqnarray*}
I_2& \leq C \Delta^{-1}\int |xF(x)|^2dx \leq C \Delta (2\pi)^{-1} \int |\gamma_1\psi_{1,\Delta}-\gamma_0\psi_{0,\Delta}| ^2.
 \end{eqnarray*}
Hence, using the factorization \eqref{decompo} we can split  $I_2\leq \pi ^{-1} C \Delta (I_{2,1}+I_{2,2})$
with 
$$\begin{cases}
   I_{2,1}=\int |\gamma_1-\gamma_0|^2,\\
I_{2,2}=\int |\gamma_0(\psi_{1,\Delta}-\psi_{0,\Delta})|^2.
  \end{cases}$$
 Using the definition of $\gamma_1$, we compute 
\begin{eqnarray}
I_{2,1}&=&c^2h_n^{2\beta+2}\int |ix_0K^*(h_n u)+h_nK^{*\prime}(h_n u)|^2du \nonumber\\
&\leq & 2c^2h_n^{2\beta+1}\left(x_0^2\int |K^*|^2+h_n^2\int |K^{*\prime}|^2\right)\nonumber \\
&\leq & 4\pi c^2h_n^{2\beta+1}\left(x_0^2\int |K|^2+h_n^2\int |xK(x)|^2\right)
\leq C h_n^{2\beta+1}. \label{I21}
\end{eqnarray}
Now, in order to deal with $I_{2,2}$, we use the previous bound \eqref{bornepsi} on  $|\psi_{1,\Delta}-\psi_{0,\Delta}|$ 
\begin{eqnarray}
I_{2,2}&\leq &C c^2\Delta^2  h_n^{2\beta+2}\int |\gamma_0(u)K^*(h_nu)|^2du 
\nonumber\\
&\leq &C c^2\Delta^2  h_n^{2\beta+2}\|K^*\|_\infty \|\gamma_0\|_2^2
\leq  C  h_n^{2\beta+1}\label{I22}
\end{eqnarray}
since $\Delta$ is bounded.

Finally, by gathering  \eqref{I1}, \eqref{I21} and \eqref{I22}, we get
\begin{eqnarray*}
 \chi^2(P_1,P_0)& \leq & C\Delta h_n^{2\beta+1}=O( n^{-1}).
 \end{eqnarray*}
This ends the proof of Theorem \ref{lower}. $\Box$

\subsection{Proof of Theorem \ref{theomet2}}

The goal is to bound $\mathbb{E}[|g(x_0)- \hat{g}_{\hat{h}}(x_0)|^2] $. To do this, we fix $ h \in H $.
We write
\begin{eqnarray*}
|g(x_0)- \hat{g}_{\hat{h}}(x_0)| \leq |\hat{g}_{\hat{h}}(x_0)- \hat{g}_{h,\hat{h}}(x_0)|+|\hat{g}_{h,\hat{h}}(x_0)- \hat{g}_{h}(x_0)|+|\hat{g}_{h}(x_0)- g(x_0)|.
\end{eqnarray*}
So we have
\begin{eqnarray*}
|g(x_0)- \hat{g}_{\hat{h}}(x_0)|^2 \leq 3|\hat{g}_{\hat{h}}(x_0)- \hat{g}_{h,\hat{h}}(x_0)|^2+ 3|\hat{g}_{h,\hat{h}}(x_0)- \hat{g}_{h}(x_0)|^2+ 3|\hat{g}_{h}(x_0)- g(x_0)|^2.
\end{eqnarray*}
Define $B:=|\hat{g}_{\hat{h}}(x_0)- \hat{g}_{h,\hat{h}}(x_0)|^2$ and $ C:= |\hat{g}_{h,\hat{h}}(x_0)- \hat{g}_{h}(x_0)|^2$.
\\
We have $A(h) \geq |\hat{g}_{\hat{h}}(x_0)- \hat{g}_{h,\hat{h}}(x_0)|^2 -V(\hat{h})\geq B -V(\hat{h})$. So $B\leq A(h) + V(\hat{h})$.\\
Moreover, $A(\hat{h})\geq |\hat{g}_{h,\hat{h}}(x_0)- \hat{g}_{h}(x_0)|^2 -V(h) \geq C-V(h)$. So $C \leq A(\hat{h}) + V(h)$.\\
Therefore,
\begin{eqnarray*}
|g(x_0)- \hat{g}_{\hat{h}}(x_0)|^2 \leq 3(A(h)+V(\hat{h}))+ 3(A(\hat{h})+ V(h))+ 3|\hat{g}_{h}(x_0)- g(x_0)|^2.
\end{eqnarray*}
Now, by definition of $\hat{h}$, $ A(\hat{h}) + V(\hat{h}) \leq A(h) + V(h)$.
This allows us to write
\begin{eqnarray*}
|g(x_0)- \hat{g}_{\hat{h}}(x_0)|^2 \leq 6A(h)+ 6 V(h)+ 3|\hat{g}_{h}(x_0)- g(x_0)|^2.
\end{eqnarray*}
Let us denote  $b_h(x_0) = \mathbb{E}[\hat{g}_{h}(x_0)]- g(x_0)$ and $b_{h,2}(x_0)= \mathbb{E}[\hat{g}_{h}(x_0)]-K_h\star g(x_0)$ (these are the same notation as in Lemma~\ref{monlem1}, but with subscript $h$). Thus
\begin{eqnarray*}
\mathbb{E}[|g(x_0)- \hat{g}_{\hat{h}}(x_0)|^2]  &\leq& 6\mathbb{E}[A(h)]+6V(h)+ 3b_h^2(x_0)+ 3\Var(\hat{g}_{h}(x_0)) \\ 
&\leq& 6\mathbb{E}[A(h)]+ 3b_h^2(x_0)+ C_2 V(h).
\end{eqnarray*}
It remains to bound $\mathbb{E}[A(h)]$. Let us denote by $g_{h,h'}= \mathbb{E}[\hat{g}_{h,h'}]$ and $g_{h}= \mathbb{E}[\hat{g}_{h}] $.
We write
\begin{equation} \label{decomp} \hat{g}_{h,h'} - \hat{g}_{h'} = \hat{g}_{h,h'} - g_{h,h'} - \hat{g}_{h'} + g_{h'} + g_{h,h'} -  g_{h'},\end{equation}
and we study the last term of the above decomposition.
We have
\begin{eqnarray*} \nonumber
|g_{h,h'}(x_0) -  g_{h'}(x_0)| & = & |\mathbb{E}[\hat{g}_{h,h'}(x_0) -  \hat{g}_{h'}(x_0)]| \\ \nonumber &=& |\mathbb{E}[K_{h'}\star\hat{g}_{h}(x_0) -  \hat{g}_{h'}(x_0)]| \\ \label{asma} &=& | K_{h'}\star\mathbb{E}[\hat{g}_{h}(x_0)- g(x_0)] +  K_{h'}\star g(x_0) - \mathbb{E}[\hat{g}_{h'}(x_0)]|.
\end{eqnarray*} 
This can be written:
\begin{eqnarray*}
|g_{h,h'}(x_0) -  g_{h'}(x_0)| & = &  | K_{h'}\star b_h(x_0) +  b_{h,2}(x_0)|\\
&\leq & {\left|\int K\left( \frac{x_0-u}{h'} \right)
 b_h(u) \frac{du}{h'}\right|} + | b_{h,2}(x_0)|.
\end{eqnarray*}
Now $| b_{h,2}(x_0)|\leq | b_{h}(x_0)|\leq \|b_h\|_\infty$ so that 
\begin{eqnarray}\nonumber
|g_{h,h'}(x_0) -  g_{h'}(x_0)|^2 
&\leq & 2\|b_h\|_\infty^2\left( \int |K(v)|dv \right)^2 +  2| b_{h,2}(x_0)|^2 \\ &\leq& 
2(\| K \|_1^2+1)\|b_h\|_\infty^2. \label{biais}
\end{eqnarray}
Then by inserting (\ref{biais}) in decomposition (\ref{decomp}), we find:
\begin{eqnarray}\nonumber
A(h) &=& \sup_{h'} {\lbrace |\hat{g}_{h,h'}(x_0) - \hat{g}_{h'}(x_0) |^2 - V(h')\rbrace}_+ \\\nonumber &\leq& 3\sup_{h'} {\lbrace  |\hat{g}_{h,h'}(x_0) - {g}_{h,h'}(x_0) |^2 - V(h')/6\rbrace}_+ \\ \label{term3} && + 3\sup_{h'} {\lbrace |\hat{g}_{h'}(x_0) - {g}_{h'}(x_0) |^2 - V(h')/6\rbrace}_+ +  6(\| K \|_1^2+1)\|b_h\|_\infty^2 .
\end{eqnarray}
We can prove the following concentration result:
\begin{prop}\label{proph}
Assume that $g$ satisfies G1, G2, G3(5) , $K$ satisfies K1,  $M=O( (n\Delta)^{1/3})$ and 
take $c$ in (\ref{cprime}) such that $c\geq 16\max(1,\|K\|_\infty).$
Then \begin{eqnarray}\label{berstein1}
\mathbb{E}\left[ \sup_{h'} {\lbrace |\hat{g}_{h'}(x_0) - {g}_{h'}(x_0) |^2 - V(h')/6\rbrace}_+ \right] =O\left( \frac{\log(n\Delta)}{n\Delta}\right)  \\
\label{berstein2}
\mathbb{E}\left[ \sup_{h'} {\lbrace |\hat{g}_{h,h'}(x_0) - {g}_{h,h'}(x_0) |^2 - V(h')/6\rbrace}_+ \right] =O \left( \frac{\log(n\Delta)}{n\Delta} \right). \end{eqnarray}
\end{prop}
Inequalities (\ref{berstein1}) et (\ref{berstein2}) together with (\ref{term3}) imply 
\begin{eqnarray*}
\mathbb{E}[|g(x_0)-\hat{g}_{\hat{h}}(x_0)|^2] \leq C_1 \|b_h\|_\infty^2 + C_2 V(h) + C_3   \frac{\log(n\Delta)}{n\Delta} .
\end{eqnarray*}
This completes the proof of Theorem \ref{theomet2}. $\Box$\\

\subsection{Proof of Theorem \ref{theoEst}. }

In all this proof, we shall use the following notation: 
$$\hat{\theta}_{\Delta}(u)=\frac1n\sum_{k=1}^{n} Z_k^{\Delta}e^{i Z_k^{\Delta}u} , \quad\hat{\eta}_{\Delta}(u)=\frac1n\sum_{k=1}^{n} (Z_k^{\Delta})^{2}e^{i Z_k^{\Delta}u},$$
and $\theta_{\Delta}(u)=\E\hat{\theta}_{\Delta}(u)$,  $\eta_{\Delta}(u)=\E\hat{\eta}_{\Delta}(u).$
We also denote $f(x)=xg(x)$, so that $f^{*}(u)=i(g^{*})'(u)$ is estimated by $\hat f_{h_1}^{*}=\hat{\eta}_{\Delta}(u)K^{*}(uh_{1})$. 
Now, let $$\Omega=\{\|g^{*}-\hat g_{h_2}^{*}\|_{2}\leq \|g^{*}\|_{2}(1-1/\sqrt2) \quad\text{ and }\quad
\|f^{*}-\hat f_{h_1}^{*}\|_{1}\leq \|f^{*}\|_{1}/2\}.$$ 
The proof is decomposed in three steps. First we shall prove that the inequality is true on $\Omega$, i.e. 
$$\mathbb{E}[|g(x_0)- \hat{g}_{\hat{h}}(x_0)|^2\1_{\Omega}] \leq C\left\lbrace \inf_{h \in H}\left\lbrace \|g- \mathbb{E}[\hat{g}_{h}] \|_\infty^2 
+ \E(V(h))\right\rbrace  + \frac{\log(n\Delta)}{n\Delta}\right\rbrace.$$
The second step is to show the rough upper bound
\begin{eqnarray*}
\mathbb{E}[|g(x_0)- \hat{g}_{\hat{h}}(x_0)|^4] \leq C(n\Delta)^{2/3}.
\end{eqnarray*}
Finally we will show that $\pr(\Omega^{c})\leq C(n\Delta)^{-8/3}$. Consequently
\begin{eqnarray*}
\mathbb{E}[|g(x_0)- \hat{g}_{\hat{h}}(x_0)|^2\1_{\Omega^{c}}] \leq \sqrt{\mathbb{E}[|g(x_0)- \hat{g}_{\hat{h}}(x_0)|^4] \pr(\Omega^{c})}\leq C(n\Delta)^{-1}
\end{eqnarray*}
and the theorem is proved.
\medskip

$\bullet$ \emph{First step: }

Following the proof of Theorem~\ref{theomet2}, we can obtain 
\begin{eqnarray*}
\E\left[|g(x_0)- \hat{g}_{\hat{h}}(x_0)|^2\1_{\Omega}\right] \leq 6 \E[A(h)\1_{\Omega}]+ 3b_h^2(x_0)+ C_2 \E(V(h)).
\end{eqnarray*}
Using the definition of $A(h)$, it is 
then sufficient to prove 
\begin{eqnarray}
\mathbb{E}\left[ \sup_{h'} {\lbrace |\hat{g}_{h'}(x_0) - {g}_{h'}(x_0) |^2 - V(h')/6\rbrace}_+ \1_{\Omega}\right] =O\left( \frac{\log(n\Delta)}{n\Delta}\right)  \\
\label{inegA2}\mathbb{E}\left[ \sup_{h'} {\lbrace |\hat{g}_{h,h'}(x_0) - {g}_{h,h'}(x_0) |^2 - V(h')/6\rbrace}_+ \1_{\Omega}\right] =O \left( \frac{\log(n\Delta)}{n\Delta} \right)
\end{eqnarray}
to obtain the result. 
Now, let us remark that on $\Omega$
$$\frac12\|g^{*}\|_{2}^{2} \leq\|\hat g_{h_{2}}^{*} \|_{2}^{2}\quad\text{ and }\quad
\frac12\|f^{*}\|_{1} \leq\|\hat f_{h_{1}}^{*} \|_{1}$$
with $\|f^{*}\|_{1}=\|(g^*)'\|_1$,
so that $$C_{0}\geq \frac{c/2}{2\pi} \|K\|^2\left(\|(g^*)'\|_1+\|g^*\|_2^2\right).$$
Then, using Proposition~\ref{proph}, since $c/2\geq16\max(1,\|K\|_{\infty})$,
\begin{eqnarray*}
&&\mathbb{E}\left[ \sup_{h'} {\lbrace |\hat{g}_{h'}(x_0) - {g}_{h'}(x_0) |^2 - V(h')/6\rbrace}_+ \1_{\Omega}\right] \\
&&\leq \mathbb{E}\left[ \sup_{h'} {\lbrace |\hat{g}_{h'}(x_0) - {g}_{h'}(x_0) |^2 - \frac16
\frac{c/2}{2\pi} \|K\|^2\left(\|(\tilde{g^*})'\|_1+\|\tilde{g^*}\|_2^2\right)\frac{\log(n\Delta)}{n\Delta}\rbrace}_+ \right] \\
&&=O\left( \frac{\log(n\Delta)}{n\Delta}\right)  
\end{eqnarray*}
and 
we prove \eqref{inegA2}
in the same way.

$\bullet$ \emph{Second step: }

First, using Lemma 3.1, $|{g}_{\hat h}(x_0)-g(x_{0})|^{2}\leq \sup_{h\in H}\left(c_{1} h ^{2}+c_{1}'\Delta^{2}\right)\leq C$. Then the bias term is uniformly bounded.
Let us now study the variance term.  We can write
\begin{eqnarray*}
\hat{g}_{h}(x_0)
&=& \frac1{2\pi} \int e^{-i x_0u}K^*(u h)\frac{1}{\Delta} \hat{\theta}_{\Delta}(u)du   \\
\end{eqnarray*}
and, since all $h\in H$ is larger than $1/M$, 
\begin{eqnarray*}
|\hat{g}_{\hat h}(x_0)-g_{\hat h}(x_{0})|
&\leq&  \frac1{2\pi}\sup_{h\in H}\int |K^*(u h)|\left|\frac{\hat{\theta}_{\Delta}(u)-\theta_{\Delta}(u)}{\Delta}\right| du \\
&\leq& \frac M{2\pi}\sum_{h\in H}\int |K^*(u )|\left|\frac{\hat{\theta}_{\Delta}(u/h)-\theta_{\Delta}(u/h)}{\Delta}\right| du.
\end{eqnarray*}
With a convex inequality
\begin{eqnarray*}
|\hat{g}_{\hat h}(x_0)-g_{\hat h}(x_{0})|^{4}&\leq& \frac{M^{7}}{(2\pi)^{4}}\sum_{h\in H}\left(\int |K^*(u )|\left|\frac{\hat{\theta}_{\Delta}(u/h)-\theta_{\Delta}(u/h)}{\Delta}\right| du\right)^{4}
\end{eqnarray*}
Next, we use the following inequality (obtained with two uses of the Schwarz inequality):
\begin{eqnarray*}
&&\E\left[(\int\phi(u)du)^{4}\right]=\iiiint\E\left[\phi(u_{1})\dots\phi(u_{4})\right]du_1\dots du_{4}\\
&&\leq \iiiint\E^{1/4}\left[\phi(u_{1})^{4}\right]\dots\E^{1/4}\left[\phi(u_{4})^{4}\right]du_1\dots du_{4}
=\left(\int\E^{1/4}\left[\phi(u)^{4}\right]du\right)^{4}.
\end{eqnarray*}
Thus,   
\begin{eqnarray*}
\E\left[|\hat{g}_{\hat h}(x_0)-g_{\hat  h}(x_{0})|^{4}\right]&\leq&
\frac{M^{7}}{(2\pi)^{4}}\sum_{h\in H}
\left(\int |K^*(u)|
\E^{1/4}\left[\left|\frac{\hat{\theta}_{\Delta}(u/h)-\theta_{\Delta}(u/h)}{\Delta}\right| ^{4}\right]du\right)^{4}\\
\end{eqnarray*}
But, according to Proposition 2.3 in~\cite{CG}, under $G3(2p)$, for  $p\geq 1$,
$\Delta^{-2p}\E\left|\hat{\theta}_{\Delta}(v)-\theta_{\Delta}(v)\right|^{2p}\leq C (n\Delta)^{-p}.$
Hence, under G3(4), 
\begin{eqnarray*}
\E|\hat{g}_{\hat h}(x_0)-g_{\hat h}(x_{0})|^{4}&\leq&   C M^{7}\sum_{h\in H}
\left(\int |K^*(u)|(n\Delta)^{-1/2}du\right)^{4}\\
&\leq&   C\|K^{*}\|_1^{4} M^{8}(n\Delta)^{-2}\leq  C\|K^{*}\|_1^{4} (n\Delta)^{2/3}.
\end{eqnarray*}

$\bullet$ \emph{Third step: }

\begin{eqnarray*}
\pr(\Omega^{c})&=&\pr(\|g^{*}-\hat g_{h_2}^{*}\|_{2}> \|g^{*}\|_{2}(1-1/\sqrt2)\text{ or }\|f^{*}-\hat f_{h_1}^{*}\|_{1}> \|f^{*}\|_{1}/2)\\
&\leq&(\|g^{*}\|_{2}(1-1/\sqrt2) )^{-8}\E\|\hat g_{h_2}^{*}-g^{*}\|_{2}^{8}+ (\|f^{*}\|_{1}/2)^{-16}\E\|\hat f_{h_1}^{*}-f^{*}\|_{1}^{16}\\
&\leq&C\left(\E\|\hat g_{h_2}^{*}-g_{h_2}^{*}\|_{2}^{8}+\E\|g_{h_2}^{*}-g^{*}\|_{2}^{8}+\E\|\hat f_{h_1}^{*}-f_{h_1}^{*}\|_{1}^{16}+\E\|f_{h_1}^{*}-f^{*}\|_{1}^{16}\right).
\end{eqnarray*}
Thus we have four terms to upperbound.
\begin{description}
\item[First term]
Since $\hat{g}_{ h_2}^{*}(u)=K_{0}^*(u h_2)\hat{\theta}_{\Delta}(u)/\Delta$,
\begin{eqnarray*}
\|\hat{g}_{ h_2}^{*}-g_{h_2}^{*}\|_{2}^{2}
&=&\frac1{h_2}\int |K_{0}^*(u)|^{2}\left|\frac{\hat{\theta}_{\Delta}(u/h_2)-\theta_{\Delta}(u/h_2)}{\Delta}\right|^{2}du.
\end{eqnarray*}
Then, under $G3(8)$, 
\begin{eqnarray*}
\E\|\hat{g}_{ h_2}^{*}-g_{h_2}^{*}\|^{8}
&\leq&\frac1{h_2^4}\left(\int\E^{1/4}\left[|K_{0}^*(u)|^{8}\left|\frac{\hat{\theta}_{\Delta}(u/h_2)-\theta_{\Delta}(u/h_2)}{\Delta}\right|^{8}\right]du\right)^{4}\\
&\leq&\frac1{h_2^4}\left(\int |K_{0}^*(u)|^{2}(n\Delta)^{-1}du\right)^{4}
\leq \|K_{0}^{*}\|_{2}^{8} M^{4}(n\Delta)^{-4}
\leq16 (n\Delta)^{-8/3}.
\end{eqnarray*}

\item[Second term]
Since $g_{h_2}^{*}=K_{0}^*(u h_2)g^{*}(u)\psi_{\Delta}(u)$, we can decompose the bias into
$$g^{*}(u)-g_{h_2}^{*}(u)=g^{*}(u)(1-K_{0}^*(u h_2))+g^{*}(u)K_{0}^*(u h_2)(1-\psi_{\Delta}(u))=b_1+b_{2}$$
Using that  $\int |g^{*}(u)|^{2}u^{2}du<\infty$, 
\begin{eqnarray*}
\|b_{1}\|^{2}&=&\int |g^{*}(u)(1-K_{0}^*(u h_2))|^{2}du=\int |g^{*}(u)|^{2}\1_{|uh_2|>1}du\\
&\leq &\int |g^{*}(u)|^{2}|uh_2|^{2}du\leq C h_2^{2}.
\end{eqnarray*}
On the other hand, using that $|1-\psi_{\Delta}(u)|\leq |u|\Delta \|g\|_{1}$ (see Proposition 2.3 in~\cite{CG})
\begin{eqnarray*}
\|b_{2}\|^{2}&=&\int |g^{*}(u)K_{0}^*(u h_2)(1-\psi_{\Delta}(u))|^{2}du\leq C\Delta^{2}\int |g^{*}(u)u|^{2}du\\
&\leq & C \Delta^{2}\leq C(n\Delta)^{-1}.
\end{eqnarray*}
Thus, taking $h_2=(n\Delta)^{-1/3}$ gives
$\|g^{*}-g_{h_2}^{*}\|^{8}\leq C h_{2}^{8}+C(n\Delta)^{-4}\leq C (n\Delta)^{-8/3}.$

\medskip
\item[Third term]
Since $\hat{f}_{h_1}^{*}(u)=K_0^*(u h_1)\hat{\eta}_{\Delta}(u)/\Delta$,  
\begin{eqnarray*}
\|\hat{f}_{ h_1}^{*}-f_{h_1}^{*}\|_{1}
&\leq&\frac1{h_1}\int |K_0^*(u)|\left|\frac{\hat{\eta}_{\Delta}(u/h_1)-\eta_{\Delta}(u/h_1)}{\Delta}\right|du \\
\end{eqnarray*}
Next, we use the following inequality
\begin{eqnarray*}
\E\left[(\int\phi(u)du)^{16}\right]\leq\left(\int\E^{1/16}\left[\phi(u)^{16}\right]du\right)^{16}.
\end{eqnarray*}
Exactly as in~\cite{CG}, using the Rosenthal inequality, we can prove
under $G3(4p)$, for  $p\geq 1$,
$\Delta^{-2p}\E\left|\hat{\eta}_{\Delta}(v)-\eta_{\Delta}(v)\right|^{2p}\leq C (n\Delta)^{-p}.$ 
Then, under $G3(32)$,
\begin{eqnarray*}
\E\|\hat{f}_{ h_1}^{*}-f_{h_1}^{*}\|_{1}^{16}
&\leq&\frac1{h_1^{16}}\left(\int\E^{1/16}\left[|K_0^*(u)|^{16}\left|\frac{\hat{\eta}_{\Delta}(u/h_1)-\eta_{\Delta}(u/h_1)}{\Delta}\right|^{16}\right]du\right)^{16}\\
&\leq&\frac1{h_1^{16}}\left(\int |K_0^*(u)|(n\Delta)^{-1/2}du\right)^{16}
\leq C \|K^*\|_1(n\Delta)^{-8/3}
\end{eqnarray*}
since $h_1=(n\Delta)^{-1/3}$.

\item[Fourth term]
Since $\eta_{\Delta}=-\psi_{\Delta}''=\Delta f^{*}\psi_{\Delta}+\Delta^{2}(g^{*})^{2}\psi_{\Delta}$, we can decompose the bias into
\begin{eqnarray*}
f^{*}(u)-f_{h_1}^{*}(u)&=&f^{*}(u)-K_0^*(u h_1) f^{*}(u)\psi_{\Delta}(u)-\Delta K_0^*(u h_1)(g^{*}(u))^{2}\psi_{\Delta}(u)\\
&=&f^{*}(u)(1-K_0^*(u h_1) )+f^{*}(u)K_0^{*}(uh_1)(1-\psi_{\Delta}(u))\\
&&-\Delta K_0^*(u h_1)(g^{*}(u))^{2}\psi_{\Delta}(u)\\
&=&b_{1}+b_{2}+b_{3}
\end{eqnarray*}
Since  $\int |f^{*}(u)|^{2}u^{2}du<\infty$, 
\begin{eqnarray*}
\|b_{1}\|_{1}&\leq&\int |f^{*}(u)(1-K_0^*(u h_1))|du=\int |f^{*}(u)|\1_{|uh_1|>1}du\\
&\leq &\left(\int |f^{*}(u)|^{2}|uh_1|^{2}du\int |uh_1|^{-2}\1_{|uh_1|>1}du\right)^{1/2}
\leq C h_1^{1/2}
\end{eqnarray*}
On the other hand, using that $|1-\psi_{\Delta}(u)|\leq |u|\Delta \|g\|_{1}$
\begin{eqnarray*}
\|b_{2}\|_{1}&\leq&\int |f^{*}(u)K_0^*(u h_1)(1-\psi_{\Delta}(u))|du\leq C\Delta \int
|f^{*}(u)uK_0^*(u h_1)|du\\
&\leq & C \Delta \left(\int |f^{*}(u)u|^{2}du\int |K_0^*(u h_1)|^{2} du\right)^{1/2}\\
&\leq & C \Delta h_{1}^{-1/2}\leq C(h_{1}n\Delta)^{-1/2}
\end{eqnarray*}
and
\begin{eqnarray*}
\|b_{3}\|_{1}&\leq&\Delta\int |K_0^*(u h_1)(g^{*}(u))^{2}\psi_{\Delta}(u)|du\\
&\leq&\Delta\int |(g^{*}(u))^{2}|du\leq  C \Delta\leq C(n\Delta)^{-1/2}
\end{eqnarray*}
Thus
$\|f^{*}-f_{h_1}^{*}\|_{1}^{16}\leq Ch_1^{8}+C(h_{1}n\Delta)^{-8}+C(n\Delta)^{-8}\leq C(n\Delta)^{-8/3}.$

\end{description}

This completes the proof of Theorem \ref{theoEst}. $\Box$\\

\subsection{Proof of Proposition \ref{proph}. }

Note that
\begin{eqnarray}\label{somme1}
&&\hat{g}_{h'}(x_0) - {g}_{h'}(x_0)  = \frac{1}{n} \sum_{k=1}^{n} 
\left[\frac{{Z_k}^\Delta}{\Delta} K_{h'} \left( x_0-{Z_k}^\Delta\right)- \mathbb{E}\left(\frac{{Z_k}^\Delta}{\Delta } K_{h'} \left( {x_0-{Z_k}^\Delta}\right)\right)\right].
\end{eqnarray}
In order to apply a Bernstein inequality, since the $  Z_k^\Delta $'s are not bounded, we truncate these variables and consider the following decomposition:
\begin{eqnarray*}
{\lbrace|{Z_k}^\Delta| \leq \mu_n\rbrace} \mbox{ and   } {\lbrace|{Z_k}^\Delta| > \mu_n\rbrace}
\end{eqnarray*}
where \begin{equation} \label{muntilde}
\mu_n  =\mu_n(h')=\frac{ \|K\|_2^2 (\|(g^*)'\|_1+ \|g^{*}\|_2^2 )}{2\pi\| K \|_{\infty} \sqrt{V(h')/6}} .
\end{equation}
We then decompose (\ref{somme1}) as follows
\begin{eqnarray*}
\hat{g}_{h'}(x_0) - {g}_{h'}(x_0)  &=& \frac{1}{n} \sum_{k=1}^{n} W_k(h')+T_k(h')-\E\left( W_k(h')+T_k(h')\right)\\
&=& S_n(W(h'))+S_n(T(h'))
\end{eqnarray*}
where $S_n(X)$ means $(1/n) \sum_{i=1}^n [X_i-\mathbb{E}(X_i)]$ and 
\begin{eqnarray}
 \label{Wh}
W_k (h)= \frac{{Z_k}^\Delta}{\Delta} K_{h} \left( {x_0-{Z_k}^\Delta}\right)\1_{\lbrace|{Z_k}^\Delta| \leq \mu_n(h) \rbrace}\\
\label{Th}
T_k (h)= \frac{{Z_k}^\Delta}{\Delta } K_{h} \left( {x_0-{Z_k}^\Delta}\right)\1_{\lbrace|{Z_k}^\Delta| > \mu_n(h) \rbrace}.
\end{eqnarray}
Thus 
\begin{eqnarray*}
  && \mathbb{E}\left[ \sup_{h'} {\lbrace |\hat{g}_{h'}(x_0) - {g}_{h'}(x_0) |^2 - V(h')/6\rbrace}_+ \right] \\
\label{eq4} &\leq & 2 \sum_{h'\in H}\mathbb{E}\left[ S_n(W(h'))^2 - V(h')/12\right]_+ +  2 \sum_{h'\in H} \mathbb{E}\left[S_n(T(h'))^2  \right].
\end{eqnarray*}
Then we use the two following lemmas
\begin{lem} \label{lem1proph} 
Assume that $g$ satisfies G1, G2, $K$ satisfies K1, and $c\geq 16, M=O( (n\Delta)^{1/3})$.
Then there exists $C>0$ only depending on $K$ and $g$ such that
\begin{eqnarray*}
 \sum_{h\in H}\mathbb{E}\left[ S_n^2(W(h))-V(h)/12 \right]_+
\leq C\frac{\log(n\Delta)}{n\Delta}. 
\end{eqnarray*}
\end{lem}
 \begin{lem} \label{lem2proph}
Under assumptions K1, G3(5) and if $M=O( (n\Delta)^{1/3}),$
 \begin{eqnarray*}
 \sum_{h\in H} \mathbb{E}\left[S_n^2(T(h))  \right]
\leq C'\frac{1}{n\Delta}.
\end{eqnarray*}
\end{lem}
Lemmas \ref{lem1proph} and \ref{lem2proph}  yield
\begin{eqnarray*}
\mathbb{E}\left[ \sup_{h'} {\lbrace |\hat{g}_{h'}(x_0) - {g}_{h'}(x_0) |^2 - V(h')/6\rbrace}_+ \right] \leq C'' \left(  \frac{1}{n\Delta} + \frac{\log(n\Delta)}{n\Delta}\right) \end{eqnarray*}

Inegality~\eqref{berstein2} is obtained by following the same lines as for inequality \eqref{berstein1} with $K_h $ replaced by 
$ K_{h'}\star K_h$.
This ends the proof of Proposition \ref{proph}. $\Box$\\

\subsection{Proof of lemma \ref{lem1proph}. }

First, note that 
\begin{eqnarray*}
 \mathbb{E}\left[ S_n^2(W(h))-V(h)/12 \right]_+
&\leq &\int_0^\infty \pr( S_n^2(W(h))\geq V(h)/12 +x)dx\\
&\leq &\int_0^\infty V(h)\pr\left(| S_n(W(h))|\geq \sqrt{V(h)(1/12 +y)}\right)dy.
\end{eqnarray*}
Next, we recall the classical Bernstein inequality (see e.g. \cite{BM98} for a proof): 
\begin{lem}
Let $W_1,...,W_n$ $n$  independent and identically distributed random variables and  $S_n(W)=(1/n) \sum_{i=1}^n [W_i-\mathbb{E}(W_i)]$. Then, for $\eta>0 $,
\begin{eqnarray*} \nonumber
\mathbb{P}(|S_n(W)| \geq \eta) &\leq & 2\exp \left( \frac{-n\eta^2 /2}{\nu^2 + b\eta}\right) 
\leq 2 \max \left( \exp \left( \frac{-n\eta^2}{4\nu^2} \right), \exp \left( \frac{-n\eta}{4b}\right)   \right),
\end{eqnarray*}
where
$Var(W_1)\leq \nu^2$ and $|W_1|\leq b $.
\end{lem}

We apply this form of Bernstein inequality to $W_i(h) $ defined by (\ref{Wh}) and  $\eta=  \sqrt{ ({1}/{12} + y) V(h)}. $ 
Using Lemma \ref{monlem2} and $\Delta\leq 1$, it is easy to see that  $$\Var(W_i) \leq\nu^2  := \frac{ {\|K\|^2_2} (\|(g^*)'\|_1+\|g^{*}\|^2_2)}{2\pi\Delta h} \mbox{ and } |W_i| \leq b:= \frac{\|K\|_{\infty} \mu_n(h)}{\Delta h} . $$
We find
 \begin{eqnarray*}
 \exp \left( \frac{-n\eta^2}{4\nu^2} \right) &=& \exp \left( - \frac{\pi(1/12) V(h) n\Delta  h}{2 \|K\|_2^2 (\|(g^*)'\|_1+ \|g^{*}\|_2^2) }\right)\times \exp \left( -\frac{\pi yV(h) n\Delta  h}{2\|K\|_2^2 (\|(g^*)'\|_1+ \|g^{*}\|_2^2) }\right) \\ &=& (n\Delta)^{-c/48} \times (n\Delta)^{-c y /4 }\end{eqnarray*}
and
  \begin{eqnarray*}
 \exp \left( \frac{-n\eta}{4b} \right) 
 &\leq&  (n\Delta)^{-c/48}\times (n\Delta)^{-c\sqrt{y/192}}.
  \end{eqnarray*}
Then we deduce
\begin{eqnarray*}
 \mathbb{E}\left[ S_n^2(W(h))-V(h)/12 \right]_+
&\leq &\int_0^\infty V(h)(n\Delta)^{-c/48} \max\left((n\Delta)^{-c y /4 },(n\Delta)^{-c\sqrt{y/192}}\right)dy\\
&\leq & V(h)(n\Delta)^{-c/48}\left( \int_0^\infty(n\Delta)^{-cy/4}dy+\int_0^\infty(n\Delta)^{-c \sqrt{y/192}}dy\right)\\
&\leq & \frac{4}{c}V(h)(n\Delta)^{-c/48}\left(\frac1{\log(n\Delta)}+\frac{96}{c\log(n\Delta)^2}\right)
\end{eqnarray*}
using that $\int_0^\infty e^{-y/\lambda}=\lambda$ and $\int_0^\infty e^{-\sqrt{y}/\lambda}=2\lambda^2$.
Replacing $V(h)$ by its value, it gives
\begin{eqnarray*}
\sum_{h\in H} \mathbb{E}\left[ S_n^2(W(h))-V(h)/12 \right]_+
&\leq & \frac{4C_0}{c}(n\Delta)^{-1-c/48}\left(1+\frac{96}{c\log(n\Delta)}\right)\sum_{h\in H}\frac1h.
\end{eqnarray*}
Recall that $H= \lbrace  \frac{k}{M},  1 \leq k \leq M \rbrace $.  Then
 \begin{eqnarray*}
\sum_{h} \frac{1}{h} = \sum_{k=1}^M \frac{M}{k} \leq \log(M) M\leq \frac13 \log(n\Delta)(n\Delta)^{1/3}.
 \end{eqnarray*}
Finally
\begin{eqnarray*}
\sum_{h\in H} \mathbb{E}\left[ S_n^2(W(h))-V(h)/12 \right]_+
&\leq & \frac{4C_0}{3c}(n\Delta)^{-2/3-c/48}\left(\log(n\Delta)+\frac{96}{c}\right)\\
&\leq & \frac{4C_0}{3c}(n\Delta)^{-1}\left(\log(n\Delta)+\frac{96}{c}\right)
\end{eqnarray*}
as soon as $c\geq 16$.
This completes the proof of lemma \ref{lem1proph}. $\Box$ \\

\subsection{Proof of lemma \ref{lem2proph}. }

For a fixed bandwidth $h$ in $H$, we can establish the following bound:
 \begin{eqnarray*}
 \mathbb{E}\left[  |S_n(T(h)) |^2  \right] &=& \Var\left[ \frac{1}{n}  \sum_{k=1}^n \frac{Z_k^\Delta}{\Delta h} K\left(  \frac{x_0-Z_k^{\Delta}}{h}\right) \1_{\lbrace{|Z_k^\Delta| > \mu_n}\rbrace} \right] \\ &\leq &  \frac{1}{n} \frac{\|K\|_{\infty}^2}{(\Delta h)^2} \mathbb{E}[({Z_1^\Delta})^2 \1_{\lbrace{|Z_1^\Delta| > \mu_n}\rbrace}] \\ &\leq& \frac{1}{n\Delta} \frac{\|K\|_{\infty}^2}{ h^2} \frac{\mathbb{E}[{{|Z_1^\Delta|}}^{w+2} / \Delta]}{{\mu_n}^w} 
\end{eqnarray*}
for any $w>0$. Recall that, according to Proposition~\ref{proposable}, $\mathbb{E}[{{|Z_1^\Delta|}}^{w+2} / \Delta]$ is bounded under G3($w+2$). 
We search conditions for $\sum_{h} h^{-2} {\mu_n}^{-w} \leq \mbox{constant}$. The following equalities hold up to constants:
 \begin{eqnarray*}
 \sum_{h\in H} \frac{1}{h^2 {\mu_n}^w}  = \sum_{h} \frac{V(h)^{w/2}}{h^2 } = \frac{\log(n\Delta)^{w/2}}{(n\Delta)^{w/2}}\sum_{h} \frac{1}{h^{2+ w/2 }}.
\end{eqnarray*}
Since $h= k/M$, this provides
\begin{eqnarray*}
\sum_{h} \frac{1}{h^{2+ w/2 }} = \sum_{k=1}^{M} \left(\frac{M}{k}\right)^{2+ w/2} =  M^{2+ w/2} \sum_{k=1}^{M} {\frac{1}{k^{2+ w/2}}}  =O( M^{2+ w/2}).
\end{eqnarray*}
Finally, as $M=O( (n\Delta)^{1/3} ) $, we have
 \begin{eqnarray*}
 \sum_{h} \frac{1}{h^2 {\mu_n}^w} \leq C\frac{M^{2+ w/2} \log(n\Delta)^{w/2}}{(n\Delta)^{w/2}} \leq C \log(n\Delta)^{w/2} (n\Delta)^{\frac{1}{3}(2+\frac{w}{2}) - \frac{w}{2} }.
\end{eqnarray*}
We need that $(2+{w}/{2})\times {1}/{3}- {w}/{2} < 0$, so we need the $Z_i $ admit a moment of order $w+2 \geq 5$.
This completes the proof of lemma \ref{lem2proph}. $\Box$ \\
%

\section*{Acknowledgement}
The authors thank Fabienne Comte and Valentine Genon-Catalot for enlightening discussions and helpful advices.
\bibliographystyle{apalike}
\bibliography{biblio2}

\end{document}